# ON THE CENTRAL CRITICAL VALUE OF THE TRIPLE PRODUCT $L$-FUNCTION

S. BÖCHERER AND R. SCHULZE-PILLOT

ABSTRACT. We compute the central critical value of the triple product $L$-function associated to three cusp forms $f_1, f_2, f_3$ with trivial character for groups $\Gamma_0(N_i)$ with square free levels $N_i$ not all of which are 1 and weights $k_i$ satisfying $k_1 \geq k_2 \geq k_3$ and $k_1 < k_2 + k_3$. This generalizes work of Gross and Kudla and gives an alternative classical proof of their results in the case $N_1 = N_2 = N_3$ with $k_1 = k_2 = k_3 = 2$.

## Introduction

Starting from the work of Garrett and of Piatetskii-Shapiro and Rallis on integral representations of the triple product $L$-function associated to three elliptic cusp forms the critical values of these $L$-functions have been studied in recent years from different points of view. From the classical point of view there are the works of Garrett [9], Satoh [22], Orloff [21], from an adelic point of view the problem has been treated by Garrett and Harris [10], Harris and Kudla [12] and Gross and Kudla [11]. Of course the central critical value is of particular interest. Harris and Kudla used the Siegel-Weil theorem to show that the central critical value is a square up to certain factors (Petersson norms and factors arising at the bad and the archimedean primes); the delicate question of the computation of the factors for the bad primes was left open. In the special situation that all three cusp forms are newforms of weight 2 and for the group $\Gamma_0(N)$ with square free level $N > 1$ Gross and Kudla gave for the first time a completely explicit treatment of this $L$-function including Euler factors for the bad places; they proved the functional equation and showed that the central critical value is a square up to elementary factors (that are explicitly given).

We reconsider the central critical value from a classical point of view, dealing with the situation of three cusp forms $f_1, f_2, f_3$ of weights $k_i$ $(i = 1, \ldots, 3)$ that are newforms for groups $\Gamma_0(N_i)$ with $N = \mathrm{lcm}(N_i)$ a squarefree integer $\neq 1$. The weights $k_i$ are subject to the restriction $k_1 < k_2 + k_3$ where $k_1 \geq \max(k_2, k_3)$; the distinction whether this inequality holds or not played an important role in [11] and [12] too. We start from the simplest possible Eisenstein series $\mathbf{E}$ of weight 2 for $\Gamma_0^{(3)}(N)$ on the Siegel space $\mathbf{H_3}$ ("summation over $C \equiv 0 \bmod N$"). After applying a suitable differential operator (depending on the weights $k_i$) to $\mathbf{E}$ we proceed in a way similar to Garrett's original approach: We restrict the differentiated Eisenstein series in a first step to $\mathbf{H_1} \times \mathbf{H_2}$ and integrate against $f_1$, the resulting function on $\mathbf{H_2}$ is then restricted to the diagonal and integrated against $f_2, f_3$. The necessary modifications to Garrett's coset decompositions (that were for level 1) are not

Both authors were supported by MSRI, Berkeley (NSF-grant DMS-9022140). R. Schulze-Pillot also supported by the Deutsche Forschungsgemeinschaft and by the Max-Planck-Institut für Mathematik, Bonn.

difficult (for the first step they have already been carried out in [2]). The actual computation of the integral is elementary and needs only standard results from the theory of newforms. It yields a Dirichlet series (2.41) whose Euler product decomposition is then computed in Section 3. The cases that $p$ divides one, two or all three of the levels $N_i$ or is coprime to $N$ must all be treated separately, which makes the discussion somewhat lengthy. However, the actual computation in each of these cases is again fairly straightforward. In Section 4 we show that the Euler factors defined in Section 3 are the "right ones" by proving the functional equation. In order to exhibit the central critical value as a square (up to elementary factors) we follow a similar strategy as [12]: the Eisenstein series $\mathbf{E}$ at $s = 0$ is expressed as a linear combination of genus theta series of quaternary positive definite integral quadratic forms. At most one of these genera (depending on the levels $N_i$ and the eigenvalues of the $f_i$ under the Atkin-Lehner involutions) contributes to the integral. Eichler's correspondence between cusp forms for $\Gamma_0(N)$ and automorphic forms on definite quaternion algebras allows then to express this contribution as an (explicitly computable) square of an element of the coefficient field of the $f_i$; this element arises as value of a trilinear form on a space of automorphic forms on the quaternion algebra and may be interpreted as value of a height pairing similar to [11].

It may be of interest to compare the advantages of the different methods applied to this problem. Although the adelic method makes it easier to obtain general results, the explicit computations needed here appear to become somewhat simpler in the classical context. In particular, by making use of the theory of newforms and of orthogonality relations for the theta series involved from [2] we can use the same Eisenstein series $\mathbf{E}$ independent of the $f_i$. This is of advantage since the pullback formalism is especially simple for this type of Eisenstein series and leads to the remarkably simple computations in Sectiions 2 and 3.

Most of this article was written while both authors were guests of the MSRI in Berkeley during its special year on automorphic forms. We wish to thank the MSRI for its hospitality and financial support. R. Schulze-Pillot was also supported by Deutsche Forschungsgemeinschaft during a visit of one month at MSRI and was a guest of the Max-Planck-Institut für Mathematik in Bonn in the final stage of the preparation of this manuscript.

*Notations* We use some standard notations from the theory of modular forms, in particular, we denote by $\mathbf{H_n}$ Siegel's upper half space of degree $n$ (the subscript $n = 1$ will be omitted); for functions $f$ on $\mathbf{H}$ and $g = \begin{pmatrix} a & b \\ c & d \end{pmatrix}$ we use $(f \mid_k g)(z) = \det(g)^{\frac{k}{2}}(cz+d)^{-k}f(g<z>)$ and similarly for the action of double cosets (Hecke operators). The operators $T(p)$ and $U(p)$ however will be used in their standard normalisation. The space of cusp forms of weight $k$ for $\Gamma_0(N) = \left\{ \begin{pmatrix} a & b \\ c & d \end{pmatrix} \in Sl_2(\mathbf{Z}) \mid \mathbf{c} \equiv \mathbf{0} \mod \mathbf{N} \right\}$ will be denoted by $[\Gamma_0(N), k]_0$.

## 1. Differential Operators

We have to deal with two types of embeddings of products of upper half spaces into $\mathbf{H_3}$ namely

$$\iota_{12} : \begin{cases} \mathbf{H} \times \mathbf{H_2} & \longrightarrow & \mathbf{H_3} \\ (z, Z) & \longmapsto & \begin{pmatrix} z & 0 \\ 0 & Z \end{pmatrix} \end{cases}$$

and

$$\iota_{111} : \begin{cases} \mathbf{H^3} & \longrightarrow & \mathbf{H_3} \\ (z_1, z_2, z_3) & \longmapsto & \begin{pmatrix} z_1 & & \\ & z_2 & \\ & & z_3 \end{pmatrix} \end{cases}$$

Without any danger of confusion we may denote by the same symbols the corresponding "diagonal" embeddings of groups:

$$\iota_{12} : Sl_2 \times Sp(2) \to Sp(3) \qquad \text{and} \qquad \iota_{111} : Sl_2^3 \to Sp(3)$$

One might try to apply Ibukiyama-type differential operators [16] in the integral representation of the triple L-functions (equivariant for $Sl_2 \times Sl_2 \times Sl_2 \hookrightarrow Sp(3)$). However in the actual computation of the integral, it is more convenient to have equivariance for $Sl_2 \times Sp(2) \hookrightarrow Sp(3)$. Therefore we use Maaß-type operators (see [20]) and the holomorphic differential operators introduced in [6]; we describe these operators here only for $Sp(3)$, but of course they also make sense for $Sp(n)$.

We start from a natural number $r$ and three (even) weights $k_1, k_2, k_3$ with $k_1 = \max\{k_i\}$ and satisfying the condition

$$k_2 + k_3 - k_1 \geq r \tag{1.1}$$

Then we define nonnegative integers $a, b, \nu_2, \nu_3$ by

$$\begin{aligned} r + a &= k_2 + k_3 - k_1 \\ k_1 &= r + a + b \\ k_2 &= r + a + \nu_2 \\ k_3 &= r + a + \nu_3 \end{aligned} \tag{1.2}$$

Then we have

$$b = \nu_2 + \nu_3 \tag{1.3}$$

We use two types of differential operators on $\mathbf{H_3}$. The first one is the Maaß operator

$$\begin{aligned} \mathcal{M}_\alpha &= \det(Z - \bar{Z})^{2-\alpha} \det(\partial_{ij}) \det(Z - \bar{Z})^{\alpha-1} \\ &= \sum_{\mu=0}^{3} \frac{\varepsilon_3(\alpha)}{\varepsilon_\mu(\alpha)} \cdot \mathrm{tr}\left( \left(Z - \bar{Z}\right)^{[\mu]} \cdot (\partial_{ij})^{[\mu]} \right) \\ &= \varepsilon_3(\alpha) + \cdots + \det(Z - \bar{Z}) \cdot \det(\partial_{ij}) \end{aligned} \tag{1.4}$$

where (following [20])

$$\varepsilon_\mu(\alpha) = \begin{cases} 1 & \mu = 0 \\ \alpha \cdot (\alpha - \tfrac{1}{2}) \cdots (\alpha - \tfrac{\mu-1}{2}) & \mu > 0 \end{cases}$$

and for a matrix $A$ of size $n$ we denote by $A^{[\mu]}$ the matrix of $\mu \times \mu$-minors. We put

$$\mathcal{M}_\alpha^{[\nu]} = \mathcal{M}_{\alpha+\nu-1} \circ \cdots \circ \mathcal{M}_{\alpha+1} \circ \mathcal{M}_\alpha$$

We recall from [20] that

$$\mathcal{M}_\alpha^{[\mu]}\left(f \mid_{\alpha,\beta} g\right) = \left(\mathcal{M}_\alpha^{[\mu]}f\right) \mid_{\alpha+\mu,\beta-\mu} g \tag{1.5}$$

for all $g = \begin{pmatrix} a & b \\ c & d \end{pmatrix} \in Sp(3,\mathbf{R})$ with $(f \mid_{\alpha,\beta} g)(Z) = \det(cZ+d)^{-\alpha} \det(c\overline{Z}+d)^{-\beta} f(g<Z>)$. Here $\alpha$ and $\beta$ are arbitrary complex numbers, but it would be sufficient for us to take $\alpha = r + s$, $\beta = s$ with $s \in \mathbf{C}$.

The second type of differential operators was introduced in [6]: It maps scalar-valued functions on $\mathbf{H_3}$ to vector-valued functions , more precisely to $\mathbf{C[X_2, X_3]_b}$ -valued functions on $\mathbf{H} \times \mathbf{H_2} \hookrightarrow \mathbf{H_3}$ where $\mathbf{C[X_2, X_3]_b}$ denotes the space of homogeneous polynomials of degree b; we realize the symmetric tensor representation $\sigma_b$ of $Gl(2,\mathbf{C})$ on this space in the usual way. The operator $\mathbf{L}_\alpha^{(\mathbf{b})}$ as defined in [6] satisfies

$$\begin{aligned}
\left(\mathbf{L}_\alpha^{(\mathbf{b})}\mathbf{f}\right) \mid_{\alpha+b,\beta}^w g_1 &= \mathbf{L}_\alpha^{(\mathbf{b})}\left(\mathbf{f} \mid_{\alpha,\beta} \iota_{1,2}(\mathbf{g_1}, \mathbf{1_4})\right)(\iota_{1,2}(\mathbf{w}, \mathbf{Z})) \\
\left(\mathbf{L}_\alpha^{(\mathbf{b})}\mathbf{f}\right) \mid_{\alpha,\beta,\sigma_b}^Z g_2 &= \mathbf{L}_\alpha^{(\mathbf{b})}\left(\mathbf{f} \mid_{\alpha,\beta} \iota_{1,2}(\mathbf{1_2}, \mathbf{g_2})\right)(\iota_{1,2}(\mathbf{w}, \mathbf{Z}))
\end{aligned} \tag{1.6}$$

for all $g_1 \in Sl_2(\mathbf{R})$ and all $g_2 \in Sp(2,\mathbf{R})$, where the upper indices $Z$ and $w$ indicate which variable is relevant at the moment and

$$\left(\left(\mathbf{L}_\alpha^{(\mathbf{b})}\mathbf{f}\right) \mid_{\alpha,\beta,\sigma_b}^Z g_2\right)(\iota_{1,2}(w, Z)) =$$
$$\det(cZ + d)^{-\alpha} \det(c\overline{Z} + d)^{-\beta}\sigma_b(cZ + d)^{-1}\left(\mathbf{L}_\alpha^{(\mathbf{b})}\right)(\iota_{12}(w, g_2 < Z >))$$

This differential operator can be described explicitly as follows:

$$\mathbf{L}_\alpha^{(\mathbf{b})} = \frac{\mathbf{1}}{\alpha^{[\mathbf{b}]}}\iota^\star \sum_{\mathbf{0 \leq 2\nu \leq b}} \frac{\mathbf{1}}{\nu!(\mathbf{b-2\nu})!(\mathbf{2 - \alpha - b})^{[\nu]}} \cdot (\mathbf{D_\uparrow D_\downarrow})^\nu \cdot (\mathbf{D - D_\uparrow - D_\downarrow})^{\mathbf{b-2\nu}} \tag{1.7}$$

with $\iota^\star$ denoting the restriction to $\mathbf{H} \times \mathbf{H_2} \hookrightarrow \mathbf{H_3}$,

$$\begin{aligned}
D_\uparrow &= \partial_{11} \\
D_\downarrow &= \sum_{2 \leq i,j \leq 3} \partial_{ij} X_i X_j \\
D - D_\uparrow - D_\downarrow &= 2\left(\partial_{12} X_1 + \partial_{13} X_3\right)
\end{aligned}$$

and

$$\alpha^{[\nu]} = \frac{\Gamma(\alpha + \nu)}{\Gamma(\alpha)} = \begin{cases} 1 & \nu = 0 \\ \alpha(\alpha + 1) \dots (\alpha + \nu - 1) & \nu > 0 \end{cases}$$

We should remark here that $\mathbf{L}_\alpha^{(\mathbf{b})}$ has coefficients, which are rational functions of $\alpha$ with no poles for $\Re(s) > 0$.

We shall use the operators

$$\mathcal{D}_\alpha^{(a,b)} := \mathbf{L}_{\alpha+\mathbf{a}'}^{(\mathbf{b})} \circ \mathcal{M}_\alpha^{[\mathbf{a}']}$$

with $2a' = a$ and $\mathcal{D}_\alpha^{\star(a,b)}$ defined by

$$\mathcal{D}_\alpha^{\star(a,b)} f = \left( \mathcal{D}_\alpha^{(a,b)} f \right) \boldsymbol{\iota}_{111}(z_1, z_2, z_3)$$

Denoting by $\mathcal{D}_\alpha^{\star(a,\nu_2,\nu_3)}$ the operator which picks out of $\mathcal{D}_\alpha^{\star(a,b)}$ its $X_2^{\nu_2} X_3^{\nu_3}$- component, we get a decomposition

$$\mathcal{D}_\alpha^{\star(a,b)} f = \sum_{\nu_2+\nu_3=b} \left( \mathcal{D}_\alpha^{\star(a,\nu_2,\nu_3)} f \right) X_2^{\nu_2} X_3^{\nu_3} \tag{1.8}$$

with

$$\begin{aligned}
&\left( \mathcal{D}_\alpha^{\star(a,\nu_2,\nu_3)} f \mid_{\alpha,\beta} \iota_{111}(g_1, g_2, g_3) \right) \\
&= \left( \mathcal{D}_\alpha^{\star(a,\nu_2,\nu_3)} f \right) \mid_{\alpha+a'+b,\beta-a'}^{z_1} g_1 \mid_{\alpha+a'+\nu_2,\beta-a'}^{z_2} g_2 \mid_{\alpha+a'+\nu_3,\beta-a'}^{z_3} g_3
\end{aligned}$$

for all $(g_1, g_2, g_3) \in Sl_2(\mathbf{R})^3$.

If $f$ is a holomorphic function on $\mathbf{H_3}$, then

$$(y_1 y_2 y_3)^{-a'} \cdot \mathcal{D}_\alpha^{\star(a,\nu_2,\nu_3)}(f)$$

is a nearly holomorphic function (in the sense of Shimura) of all three variables $z_1$, $z_2$, $z_3 \in \mathbf{H}$.

To apply Shimura's results on nearly holomorphic functions, it is more convenient to use his differential operators $\delta_\alpha^\mu$, which differ (in the one-dimensional case i.e. on $\mathbf{H}$) from the Maaß operators only by a factor constant $\times y^\mu$:

$$\delta_\alpha = \frac{1}{2\pi i} \left( \frac{\alpha}{2iy} + \frac{\partial}{\partial z} \right)$$

$$\delta_\alpha^\mu = \delta_{\alpha+2\mu-2} \circ \cdots \circ \delta_\alpha$$

By elementary considerations about the degree of nearly holomorphic functions (as polynomials in $y^{-1}$) Shimura observed that nearly holomorphic functions on $\mathbf{H}$ are linear combinations of functions obtained from holomorphic functions by applying the operators $\delta_\alpha^\mu$ (at least if $\alpha$ is not in a certain finite set, for details see [24, lemma 7]. By the same kind of reasoning we get an identity

$$\begin{aligned}
&(y_1 y_2 y_3)^{-a'} \cdot \mathcal{D}_\alpha^{\star(a,\nu_2,\nu_3)} \\
&= \sum_{0 \le \mu_1, \mu_2, \mu_3 \le a'} \delta_{\alpha+a+b-2\mu_1}^{\mu_1} \delta_{\alpha+a+\nu_2-\mu_2}^{\mu_2} \delta_{\alpha+a+\nu_3-2\mu_3}^{\mu_3} \mathbf{D}_\alpha(\mathbf{a}, \nu_2, \nu_3, \mu_1, \mu_2, \mu_3)
\end{aligned} \tag{1.9}$$

We understand that in (1.9) the operator $\delta^{\mu_i}$ acts with respect to $z_i$, $i = 1, 2, 3$; moreover $\mathbf{D}_\alpha(\dots)$ is a holomorphic differential operator mapping functions on $\mathbf{H_3}$ to functions on $\mathbf{H} \times \mathbf{H} \times \mathbf{H}$. Following again the same line of reasoning as in lemma 7 [loc.cit], adapted to our situation, we easily get that the $\mathbf{D}_\alpha(\dots)$ satisfy

$$\begin{aligned}
&\mathbf{D}_\alpha(\mathbf{a}, \nu_2, \nu_3, \mu_1, \mu_2, \mu_3) \left( \mathbf{f} \mid_{\alpha,\beta} \iota_{111}(\mathbf{g_1}, \mathbf{g_2}, \mathbf{g_3}) \right) \\
&= \mathbf{D}_\alpha(\mathbf{a}, \nu_2, \nu_3, \mu_1, \mu_2, \mu_3)(\mathbf{f}) \mid_{\alpha+\mathbf{a}+\mathbf{b}-2\mu_1,\beta}^{\mathbf{z_1}} \mathbf{g_1} \mid_{\alpha+\mathbf{a}+\nu_2-2\mu_2,\beta}^{\mathbf{z_2}} \mathbf{g_2} \mid_{\alpha+\mathbf{a}+\nu_3-2\mu_3,\beta}^{\mathbf{z_3}} \mathbf{g_3}
\end{aligned} \tag{1.10}$$

for all $(g_1, g_2, g_3) \in Sl(2, \mathbf{R})^3$ and all holomorphic functions on $\mathbf{H}_3$ (and hence also for all $C^\infty$-functions). The upper indices $z_i$ on the right hand side of (1.10) indicate, on which variable $g_i$ operates.

We have to remark here that Shimura's condition "$k > 2r$" in Lemma 7 [loc.cit] is satisfied in our situation as long as $\alpha$ is non-real or

$$\begin{aligned}
\alpha + a + b &> a \\
\alpha + a + \nu_2 &> a \\
\alpha + a + \nu_3 &> a
\end{aligned}$$

However the coefficients of the $\partial_{ij}$ on both sides of (1.9) are easily seen to be *rational* functions of $\alpha$, therefore (1.9) (and subsequent equations) are true for all $\alpha \in \mathbf{C}$ as rational functions of $\alpha$.

It is crucial for us to see that in the identity (1.10) the "holomorphic part", i.e.

$$\mathbf{D}_\alpha(\mathbf{a}, \nu_2, \nu_3, \mathbf{0}, \mathbf{0}, \mathbf{0})$$

is different from zero. For this purpose we consider $(y_1 y_2 y_3)^{-a'} \cdot \mathcal{D}_\alpha^{\star(a,b)}$ as a polynomial in $\partial_{12}, \partial_{13}, \partial_{23}$. It is easy to see that this is a polynomial of total degree $3a + b$, the component of degree $3a + b$ being given by

$$(\partial_{12} \partial_{13} \partial_{23})^a \times \frac{1}{\alpha^{[b]}} \cdot \frac{1}{b!} \cdot 2^b \cdot (y_1 y_2 y_3)^{-a'} \cdot \mathcal{D}_\alpha^{\star(a,b)}$$

as a polynomial in $\partial_{12}, \partial_{13}, \partial_{23}$. It is easy to see that this is a polynomial of total degree $3a + b$, the component of degree $3a + b$ being given by

$$(\partial_{12} \partial_{13} \partial_{23})^a \times \frac{1}{\alpha^{[b]}} \cdot \frac{1}{b!} \cdot 2^b \cdot (\partial_{12} X_2 + \partial_{13} X_3)^b \tag{1.11}$$

In particular, this component is free of $y_i^{-1}$ and $\partial_{ii}$, so it can only come from the

$$\mathbf{D}_\alpha(\mathbf{a}, \nu_2, \nu_3, \mathbf{0}, \mathbf{0}, \mathbf{0})$$

with $\nu_2 + \nu_3 = b$.

Now we define a polynomial $Q_\alpha$ of the matrix variable $S = S^t = (s_{i,j})_{1 \le i,j \le 3}$ by

$$\mathbf{D}_\alpha(\mathbf{a}, \mathbf{b}, \nu_2, \nu_3, \mathbf{0}, \mathbf{0}, \mathbf{0}) \mathbf{e}^{\text{trace}(\mathbf{SZ})} = \mathbf{Q}_\alpha(\mathbf{S}) \mathbf{e}^{\mathbf{s_{11} Z_{11} + s_{22} Z_{22} + s_{33} Z_{33}}} \tag{1.12}$$

By the same kind of reasoning as in [1, Satz 15] we see that by

$$(\mathbf{x_1}, \mathbf{x_2}, \mathbf{x_3}) \longmapsto \mathbf{P_r}(\mathbf{x_1}, \mathbf{x_2}, \mathbf{x_3}) := \mathbf{Q_r} \begin{pmatrix} \mathbf{x_1^t x_1} & \mathbf{x_1^t x_2} & \mathbf{x_1^t x_3} \\ \mathbf{x_2^t x_1} & \mathbf{x_2^t x_2} & \mathbf{x_2^t x_3} \\ \mathbf{x_3^t x_1} & \mathbf{x_3^t x_2} & \mathbf{x_3^t x_3} \end{pmatrix} \tag{1.13}$$

we get a polynomial function of $(\mathbf{x_1}, \mathbf{x_2}, \mathbf{x_3}) \in (\mathbf{C^{2r}})^3$ which in each variable is a harmonic form of degree $a + b$, $a + \nu_2$, $a + \nu_3$ respectively; more precisely, $P_r$ defines a non-zero element of

$$(\mathcal{H}_{a+b}(2r) \otimes \mathcal{H}_{a+\nu_2}(2r) \otimes \mathcal{H}_{a+\nu_3}(2r))^{\mathbf{O(2r)}} \tag{1.14}$$

For our investigation of the functional equation of triple L-functions we have to modify $\mathbf{D}_\alpha(\mathbf{a}, \mathbf{b}, \boldsymbol{\nu_2}, \boldsymbol{\nu_3}, \mathbf{0}, \mathbf{0}, \mathbf{0})$ still further (we switch notation now from $\alpha$ to $r + s$) We consider the operator $\Delta = \Delta_{r,s}(a, \nu_2, \nu_3)$ given by

$$F \longmapsto (y_1 y_2 y_3)^s \mathbf{D_{r+s}}(\mathbf{a}, \boldsymbol{\nu_2}, \boldsymbol{\nu_3}, \mathbf{0}, \mathbf{0}, \mathbf{0}) \left( \mathbf{F} \times \det(\mathbf{Y})^{-s} \right) \qquad (1.15)$$

This operator (acting on functions on $\mathbf{H_3}$) is easily seen to satisfy

$$\Delta \left( F \mid_r \iota_{111}(g_1, g_2, g_3) \right) = \Delta(F) \mid^{z_1}_{r+a+b} g_1 \mid^{z_2}_{r+a+\nu_2} g_2 \mid^{z_3}_{r+a+\nu_3} g_3 \qquad (1.16)$$

for all $g_1, g_2, g_3 \in Sl_2(\mathbf{R})$. By the same kind of argument about nearly holomorphic functions as above we get

$$\Delta_{r,s}(a, \nu_2, \nu_3) = \sum_{\substack{0 \le \mu_1 \le \left[\frac{a+b}{2}\right] \\ 0 \le \mu_2 \le \left[\frac{a+\nu_2}{2}\right] \\ 0 \le \mu_3 \le \left[\frac{a+\nu_3}{2}\right]}} \delta^{\mu_1}_{r+a+b-2\mu_1} \delta^{\mu_2}_{r+a+\nu_2-2\mu_2} \delta^{\mu_3}_{r+a+\nu_3-2\mu_3} \Delta_{r,s}(a, \nu_2, \nu_3, \mu_1, \mu_2, \mu_3)$$

$$(1.17)$$

with holomorphic differential operators $\Delta_{r,s}(a, \nu_2, \nu_3, \mu_1, \mu_2, \mu_3)$ mapping functions on $\mathbf{H_3}$ to functions on $\mathbf{H} \times \mathbf{H} \times \mathbf{H}$. We should mention here that the differential operators coming up in (1.17) do not have poles as long as $r$ is positive and $s$ is non-real or $Re(s) \ge 0$. Again the "holomorphic part" $\Delta_{r,s}(a, \nu_2, \nu_3, 0, 0, 0)$ defines (as in (1.12),(1.13) an element of

$$\left( \mathcal{H}_{a+b}(2r) \otimes \mathcal{H}_{a+\nu_2}(2r) \otimes \mathcal{H}_{a+\nu_3}(2r) \right)^{\mathbf{O(2r)}}$$

This space is known to be one-dimensional: By a result of Littelmann ([19], p. 145) the decomposition of $\mathcal{H}_{a+\nu_2}(2r) \otimes \mathcal{H}_{a+\nu_3}(2r)$ is multiplicity free and contains $\mathcal{H}_{a+b}(2r)$), hence there is a unique invariant line in the threefold tensor product. There exists therefore a function $c = c_r(s)$ such that

$$\Delta_{r,s}(a, \nu_2, \nu_3, 0, 0, 0) = c_r(s) \mathbf{D_r}(\mathbf{a}, \boldsymbol{\nu_2}, \boldsymbol{\nu_3}, \mathbf{0}, \mathbf{0}, \mathbf{0}). \qquad (1.18)$$

By comparing coefficients of $(\partial_{12} \partial_{13} \partial_{23})^a \partial^{\nu_2}_{12} \partial^{\nu_3}_{13}$ on both sides of (1.17) we get

$$c_r(s) = \frac{(r+a')^{[b]}}{(r+s+a')^{[b]}} = \frac{\Gamma(r+a'+b)}{\Gamma(r+a')} \cdot \frac{\Gamma(r+a'+s)}{\Gamma(r+a'+s+b)} \qquad (1.19)$$

## 2.Unfolding the integral

For a squarefree number $N > 0$ and three cuspforms

$$f = \sum a_f(n) e^{2\pi i n z} \in [\Gamma_0(N), k_1]_0$$
$$\phi = \sum a_\phi(n) e^{2\pi i n z} \in [\Gamma_0(N), k_2]_0$$
$$\psi = \sum a_\psi(n) e^{2\pi i n z} \in [\Gamma_0(N), k_3]_0$$

with $k_1, k_2, k_3$ as in section 1 we want to compute the threefold integral $\mathcal{A}(f, \phi, \psi, s)$, defined by

$$\int\limits_{(\Gamma_0(N)\backslash \mathbf{H})^3} \int \int \overline{f(z_1)\phi(z_2)\psi(z_3)} \left(\mathcal{D}^{\star(a,\nu_2,\nu_3)}\left(\mathbf{G}^3_{\mathbf{r,s}}\right)\right)(\iota_{111}(z_1, z_2, z_3))$$

$$\times\, y_1^{k_1+s-a'} y_2^{k_2+s-a'} y_3^{k_3+s-a'} \prod_{i=1}^3 \frac{dx_i dy_i}{y_i^2} \quad (2.1)$$

where $\mathbf{G}^3_{\mathbf{r,s}}$ is the Eisenstein series on $\mathbf{H_3}$ defined by

$$\begin{aligned}\mathbf{G}^3_{\mathbf{r,s}} &= \sum_{M \in \Gamma^3_\infty \backslash \Gamma^3_0(N)} 1 \mid_{r+s,s} M \\ &= \sum_{\begin{pmatrix} * & * \\ C & D \end{pmatrix} = M \in \Gamma^3_\infty \backslash \Gamma^3_0(N)} \det(CZ+D)^{-r-s} \det(C\bar{Z}+D)^{-s} \end{aligned} \quad (2.2)$$

In the applications we shall need modified versions of the integral (2.1); it is appropriate to describe these here: We use the well-known fact (see e.g. [25, equation (2.28)]) that holomorphic cusp forms are orthogonal to $(C^\infty\text{-})$automorphic forms in the image of the differential operators $\delta$, therefore we may replace $\left(\mathcal{D}^{\star(a,\nu_2,\nu_3)}\left(\mathbf{G}^3_{\mathbf{r,s}}\right)\right)(\iota_{111}(z_1, z_2, z_3)) \, y_1^{k_1+s-a'} y_2^{k_2+s-a'} y_3^{k_3+s-a'}$ in the integrand of (2.1) by

$$\left(\mathbf{D_{r+s}(a, \nu_2, \nu_3, 0, 0, 0)}\left(\mathbf{G}^3_{\mathbf{r,s}}\right)\right)(\iota_{111}(z_1, z_2, z_3)) \, y_1^{k_1+s} y_2^{k_2+s} y_3^{k_3+s}$$

or by

$$c_r(s)\mathbf{D_r(a, \nu_2, \nu_3, 0, 0, 0)}\left(\mathbf{E}^3_{\mathbf{r,s}}\right)(\iota_{111}(\mathbf{z_1, z_2, z_3})\mathbf{y_1^{k_1} y_2^{k_2} y_3^{k_3}} \quad (2.3)$$

where

$$\mathbf{E}^3_{\mathbf{r,s}}(\mathbf{Z}) = \det(\mathbf{Y})^{\mathbf{s}} \cdot \mathbf{G}^3_{\mathbf{r,s}} = \sum_{\mathbf{M} \in \Gamma^3_\infty \backslash \Gamma^3_0(N)} \det(\mathbf{Y})^{\mathbf{s}} \mid_{\mathbf{r}} \mathbf{M}$$

The actual computation of $\mathcal{A}(f, \phi, \psi, s)$ is however most conveniently done using the integral in the version of (2.1)

## 2.1. The first integration
To understand the integration with respect to $z_1$, it is better to consider first the integral

$$I(s) = \int\limits_{\Gamma_0(N)\backslash \mathbf{H}} \overline{f(z_1)} \left(\mathcal{D}^{(a,b)}_{r+s} \mathbf{G}^3_{\mathbf{r,s}}\right)(\iota_{1,2}(z_1, Z)) \, y_1^{k_1+s-a'} \frac{dx_1 dy_1}{y_1^2} \times \det(Y)^{s-a'} \quad (2.4)$$

with $Z = X + iY \in \mathbf{H_2}$. We recall from [2, Thm.1.1] that the double cosets

$$\Gamma_\infty \backslash \Gamma^3_0(N) / \iota_{1,2}\left(\Gamma_0(N) \times \Gamma^2_0(N)\right)$$

can be parametrized by the following set of representatives

$$\{\mathbf{g_m} \mid \mathbf{m} \in \mathbf{N} \cup \{\mathbf{0}\}, \mathbf{m} \equiv \mathbf{0} \bmod \mathbf{N}\} \tag{2.5}$$

where

$$\mathbf{g_m} = \begin{pmatrix} & \mathbf{1_3} & & \mathbf{0_3} \\ 0 & m & 0 \\ m & 0 & 0 & \mathbf{1_3} \\ 0 & 0 & 0 \end{pmatrix}$$

We split the integral (2.4) into the contributions of the double cosets (2.5):

$$I(s) = \sum_m I_m(s) \tag{2.6}$$

It is easy to see that the double coset with $m = 0$ decomposes into left cosets as follows

$$\left\{ \iota_{1,2}(\gamma, \delta) \mid \gamma \in \Gamma_\infty \backslash \Gamma_0(N), \delta \in \Gamma_\infty^2 \backslash \Gamma_0^2(N) \right\}$$

Therefore its contribution to $\mathcal{D}^{\star(a,b)} \mathbf{G}^{\mathbf{3}}_{\mathbf{r+s,s}}$ is just

$$\mathcal{D}^{\star(a,b)}_{r+s,s} \left( \sum_{\gamma, \delta} 1 \mid_{r+s,s} \iota_{1,2}(\gamma, \delta) \right)$$
$$= \sum_{\gamma, \delta} \mathcal{D}^{\star(a,b)}_{r+s}(1) \mid_{r+s+a'+b,s-a'}^{z_1} \gamma \mid_{r+s+a',s-a',\sigma_b}^Z \delta$$

It is obvious that

$$\mathcal{D}^{\star(a,b)}_{r+s}(1) = \begin{cases} 0 & b > 0 \\ \varepsilon_3(r+s)\varepsilon_3(r+s+1)\dots\varepsilon_3(r+s+a'-1) & b = 0 \end{cases}$$

Unfolding the integral defining $I_0(s)$ we easily get (by the cuspidality of $f$) that

$$I_0(s) = 0$$

(we omit the standard calculation)

For fixed $m > 0$, $m \equiv 0 \bmod N$ the left cosets are given by (see [2, Thm1.2])

$$\left\{ \mathbf{g_m} \iota_{\mathbf{1,2}}(\gamma, \mathbf{l(h)g}) \mid \gamma \in \mathbf{\Gamma_0(N)}, \mathbf{h} \in \Gamma[\mathbf{m}] \backslash \mathbf{\Gamma_0(N)}, \mathbf{g} \in \mathbf{C_{2,1}(N)} \backslash \mathbf{\Gamma_0^2(N)} \right\} \tag{2.7}$$

where

$$\Gamma[m] = \Gamma_0(N) \cap \begin{pmatrix} 0 & m^{-1} \\ -m & 0 \end{pmatrix} \Gamma_0(N) \begin{pmatrix} 0 & -m^{-1} \\ m & 0 \end{pmatrix},$$

$l = \iota_{11}$ and $C_{2,1}$ is the standard maximal parabolic subgroup of $Sp(2)$ given by

$$C_{2,1} = \left\{ \begin{pmatrix} & & & \\ & & & \\ & & & \\ 0 & 0 & 0 & \end{pmatrix} \right\}$$

with

$$C_{2,1}(N) := C_{2,1}(\mathbf{Z}) \cap \mathbf{\Gamma_0^2(N)}.$$

The summation over $\gamma$ unfolds the integral for $I_m(s)$ to

$$I_m(s) = \int_{\mathbf{H}} \overline{f(z_1)} \sum_{h,g} \mathcal{D}_{r+s}^{\star(a,b)} \left(1 \mid_{r+s,s} \mathbf{g_m}\right) \mid_{r+s+a',s-a',\sigma_b}^{Z} l(h)g \cdot y_1^{k_1+s-a'} \frac{dx_1 dy_1}{y^2} \times \det(Y)^{s-a'} \quad (2.8)$$

By lemma 4.2 of [6] and (1.4) we have

$$\begin{aligned}
&\mathcal{D}^{\star(a,b)} \left(1 \mid_{r+s,s} \mathbf{g_m}\right) \\
&= 2^{-a}(2r+2s-2)^{[a]}(r+s)^{[a']} \mathbf{L}_{\mathbf{r+s+a'}}^{\mathbf{(b)}} \left(\mathbf{1} \mid_{\mathbf{r+s+a',s-a'}} \mathbf{g_m}\right) \\
&= A(r+s,b) \cdot (1-m^2 z_1 Z^*)^{-r-s-a'-b} \cdot (1-m^2 \bar{z}_1 \bar{Z}^*)^{-s+a'} \cdot (m X_2)^b
\end{aligned} \quad (2.9)$$

where for $Z \in \mathbf{H_2}$ we denote by $Z^*$ the entry in the upper left corner of $Z$ and

$$A(s,b) = \left(\frac{-1}{2\pi i}\right)^b \cdot \frac{2^{-a}(2s-2)^{[a]}(s)^{[a']}(2s+a-2)^{[b]}}{b!(s+a'-1)^{[b]}}$$

This implies

$$I_m(s) = A(r+s,b) \sum_{g,h} \left( \int_{\mathbf{H}} \overline{f(z_1)} \left(1-m^2 z_1 Z^*\right)^{-r-s-a'-b} \left(1-m^2 \bar{z}_1 \bar{Z}^*\right)^{-s+a'} \right. \\
\left. \times y_1^{k_1+s-a'} \frac{dx_1 dy_1}{y^2} \times (m X_2)^b \right) \mid_{r+s+a',s-a',\sigma_b}^{Z} l(h)g \det(Y)^{s-a'} \quad (2.10)$$

The integral in (2.10) is exactly of the same type as in [2, (1.4),(1.5)]. Using the same notation

$$\mu(k,s) = (-1)^{\frac{k}{2}} \cdot 2^{3-k-2s} \frac{\pi}{k+s-1}$$

as in [2] we get by the same reasoning as there

$$\begin{aligned}
I_m(s) = {}&A(r+s,b) \cdot \mu(r+a+b,s-a') \cdot m^{-r-2s} \\
&\times \sum_{g} \left( f^\rho \mid_{r+b+a} \Gamma_0(N) \begin{pmatrix} 0 & -m^{-1} \\ m & 0 \end{pmatrix} \Gamma_0(N) \right) (g<Z>^*> \\
&\times \left( \sigma_b \otimes \overset{k+a}{\det} (j(g,Z)) \right)^{-1} X_2^b \left( \frac{\det \Im(g<Z>)}{\Im(g<Z>^*)} \right)^{s-a'}
\end{aligned} \quad (2.11)$$

with $j(\begin{pmatrix} a & b \\ c & d \end{pmatrix}, z) = cz+d$. This is (essentially) a vector-valued Klingen-type Eisenstein series attached to the modular form

$$f^\rho \mid_{k_1} \Gamma_0(N) \begin{pmatrix} 0 & -m^{-1} \\ m & 0 \end{pmatrix} \Gamma_0(N)$$

where

$$f^\rho(z) = \overline{f(-\bar{z})}.$$

From now on we assume that $f$ is a normalized newform (eigenform) of level $N_f \mid N$; we write $N = N_f \cdot N^f$ whenever it is convenient. The Fourier coefficients of $f$ are then totally real and we have an Euler product expansion of type

$$\sum a_f(n)n^{-s} = \left(\prod_{p \mid N_f} \frac{1}{(1 - a_f(p)p^{-s})}\right) \cdot \left(\prod_{p \nmid N_f} \frac{1}{(1 - \alpha_p p^{-s})(1 - \alpha'_p p^{-s})}\right) \qquad (2.12)$$

Moreover for any prime $q$ dividing $N_f$ we have

$$a_f(q)^2 = q^{k_1 - 2}$$

and

$$f|_{k_1}\mathbf{V_q^N} = \mathbf{f}|_{\mathbf{k_1}}\mathbf{V_q^{N_f}} = -\mathbf{a_f(q)q^{1-\frac{k_1}{2}}f}$$

where $\mathbf{V_q^N}$ denotes the "Atkin-Lehner-involution" given by

$$\mathbf{V_q^N} = \begin{pmatrix} x & y \\ N & q \end{pmatrix}$$

with $xq - Ny = q$ and $q|x$ (for details we refer to [18]).

Actually we have to work not with the newform $f$ itself, but with $f|_{k_1}\begin{pmatrix} D^f & 0 \\ 0 & 1 \end{pmatrix}$ where $D^f$ is a fixed divisor of $N^f$

To simplify (2.11) further (for $f|_{k_1}\begin{pmatrix} D^f & 0 \\ 0 & 1 \end{pmatrix}$), we have to study

$$\begin{aligned}
&\sum_{m \equiv 0(N)} f\mid_{k_1}\begin{pmatrix} D^f & 0 \\ 0 & 1 \end{pmatrix}\mid_{k_1}\Gamma_0(N)\begin{pmatrix} 0 & -m^{-1} \\ m & 0 \end{pmatrix}\Gamma_0(N) \cdot m^{-s} \\
&= N^{-s}\sum_{m'=1}^{\infty} f\mid_{k_1}\begin{pmatrix} D^f & 0 \\ 0 & 1 \end{pmatrix}\mid_{k_1}\Gamma_0(N)\begin{pmatrix} 1 & 0 \\ 0 & m'^2 N \end{pmatrix}\Gamma_0(N)\mid_{k_1}\begin{pmatrix} 0 & -1 \\ N & 0 \end{pmatrix} \cdot m'^{-s}
\end{aligned} \qquad (2.13)$$

Now we are essentially in a "local" situation, because we may decompose the "Fricke involution" into Atkin-Lehner involutions:

$$\begin{pmatrix} 0 & -1 \\ N & 0 \end{pmatrix} = \gamma \circ \prod_{q \mid N} \mathbf{V_q^N}$$

with $\gamma \in \Gamma_0(N)$.

We use the following formal identities:

$p \nmid N$ :
$$\sum_{l=0}^{\infty} f \mid \Gamma_0(N)\begin{pmatrix} 1 & 0 \\ 0 & p^{2l} \end{pmatrix}\Gamma_0(N)X^l = \frac{(1-X)(1-p^2 X^2)}{(1-pX)(1-\alpha_p^2 p^{-k_1+2}X)(1-\alpha_p'^2 p^{-k_1+2}X)} \qquad (2.14)$$

*Proof*: Standard

$p|N_f$ :

$$\sum_{l=0}^{\infty} f \mid \Gamma_0(N) \begin{pmatrix} 1 & 0 \\ 0 & p^{2l+1} \end{pmatrix} \Gamma_0(N) \circ \mathbf{V_p^N} \cdot \mathbf{X^l} = \mathbf{p^{1-\frac{k_1}{2}}} \cdot \frac{\mathbf{a_f(p)}}{\mathbf{1-X}} \cdot \mathbf{f}|_{\mathbf{k_1}} \mathbf{V_p^{N_f}}$$

$$= \frac{-1}{1-X} f$$

(2.15)

*Proof*: Standard, using $a_f(p^n) = a_f(p)^n$ and $a_f(p^2) = p^{k_1-2}$

$p|N^f$ :

$$\sum_{l=0}^{\infty} f \mid \Gamma_0(N) \begin{pmatrix} 1 & 0 \\ 0 & p^{2l+1} \end{pmatrix} \Gamma_0(N) \circ \mathbf{V_p^N} \times \mathbf{X^l}$$

$$= p^{1-\frac{k_1}{2}} \cdot \frac{(f|U(p) - p^{\frac{k}{2}} \cdot f|_{k_1} \begin{pmatrix} p & 0 \\ 0 & 1 \end{pmatrix} \cdot X)}{(1-\alpha_p^2 p^{-k_1+2}X)(1-{\alpha'_p}^2 p^{-k_1+2}X)}|_{k_1} \mathbf{V_p^N}$$

(2.16)

*Proof:* Standard.
Using

$$f|T(p) = f|U(p) + p^{\frac{k_1}{2}-1} \cdot f|_{k_1} \begin{pmatrix} p & 0 \\ 0 & 1 \end{pmatrix}$$

and

$$f|_{k_1} \begin{pmatrix} p & 0 \\ 0 & 1 \end{pmatrix}|_{k_1} \mathbf{V_p^N} = \mathbf{f}$$

we get for (2.16)

$$\frac{p^{1-\frac{k_1}{2}} a_f(p) f|_{k_1} \begin{pmatrix} p & 0 \\ 0 & 1 \end{pmatrix} - (1+pX)f}{(1-\alpha_p^2 p^{-k_1+2}X)(1-{\alpha'_p}^2 p^{-k_1+2}X)}$$

(2.17)

In quite the same way we get (still for the case $p|N^f$):

$$\sum_{l=0}^{\infty} f|_{k_1} \begin{pmatrix} p & 0 \\ 0 & 1 \end{pmatrix} | \Gamma_0(N) \begin{pmatrix} 1 & 0 \\ 0 & p^{2l+1} \end{pmatrix} \Gamma_0(N)|_{k_1} \mathbf{V_p^N} \cdot \mathbf{X^l}$$

$$= \frac{p(1+pX) \cdot f - a_f(p)p^{-\frac{k_1}{2}+2}X \cdot f \mid_{k_1} \begin{pmatrix} p & 0 \\ 0 & 1 \end{pmatrix}}{(1 - \alpha_p^2 p^{-k_1+2}X)(1 - {\alpha_p'}^2 p^{-k_1+2}X)}|_{k_1} \mathbf{V_p^N} \qquad (2.18)$$

$$= \frac{p(1+pX)f|_{k_1} \begin{pmatrix} p & 0 \\ 0 & 1 \end{pmatrix} - a_f(p)p^{2-\frac{k_1}{2}}f \cdot X}{(1 - \alpha_p^2 p^{-k_1+2}X)(1 - {\alpha_p'}^2 p^{-k_1+2}X)}$$

The usual procedure ( $X \longmapsto p^{-s}$ )yields for (2.13)

$$\sum_{m \equiv 0(N)} f \mid_{k_1} \begin{pmatrix} D^f & 0 \\ 0 & 1 \end{pmatrix}|_{k_1} \Gamma_0(N) \begin{pmatrix} 0 & -m^{-1} \\ m & 0 \end{pmatrix} \Gamma_0(N) \circ \begin{pmatrix} 0 & -1 \\ N & 0 \end{pmatrix} \cdot m^{-s}$$

$$= N^{-s} \prod_{p|N} \frac{(1-p^{-s})(1-p^{2-2s})}{1-p^{-s+1}} \prod_{p|N_f} \frac{1}{(1 - \alpha_p^2 p^{-s-k_1+2})(1 - \alpha_p^2 p^{-s-k_1+2})} \prod_{p|N_f} \frac{-1}{(1-p^{-k_1-s})} \times \tilde{f}_s \qquad (2.19)$$

with

$$\tilde{f}_s = \sum_{d|N^f} \alpha(d, D, s) f|_{k_1} \begin{pmatrix} d & 0 \\ 0 & 1 \end{pmatrix} \qquad (2.20)$$

where

$$\alpha(d, D, s) = \prod_{p|N^f} \alpha(d_p, D_p, s) \qquad (2.21)$$

is a multiplicative function given by (2.14)-(2.18). Here we denote by $t_p$ the $p$-part of a positive rational number $t$. In the sequel we write $\alpha_p(d, D, s)$ instead of $\alpha(d_p, D_p, s)$.

## 2.2 Second Unfolding

To continue the computation of the integral (2.1) we first need to find a good parametrization of $C_{2,1}(N)\backslash\Gamma_0^2(N)$; we shall follow [22] (with the modifications necessary for level $N > 1$). We first remark that two elements of $Sp_2(\mathbf{Z})$ are equivalent modulo $C_{2,1}(\mathbf{Z})$ iff their last rows are equal up to sign (the same is true for $\Gamma_0(N)$ and $C_{2,1}(N)$).

For $C_{2,1}(\mathbf{Z})\backslash\mathbf{Sp_2}(\mathbf{Z})$ the parametrization given in [22] is as follows:

$$\{\iota_{1,1}(\mathbf{1_2}, h) \mid h \in Sl_2(\mathbf{Z})_\infty\backslash\mathbf{Sl_2}(\mathbf{Z})\} \qquad (2.22)$$

$$\left\{d(\mathbf{J})\iota_{1,1}(\mathbf{h}, \mathbf{1_2}) \mid \mathbf{h} \in \mathbf{Sl_2}(\mathbf{Z})_\infty^+\backslash\mathbf{Sl_2}(\mathbf{Z})\right\} \qquad (2.23)$$

$$\bigcup_{\substack{u, v \in \mathbf{N} \\ u, v \text{ coprime}}} \left\{d(M) \circ \iota_{1,1}(h, h') \mid h \in Sl_2(\mathbf{Z})_\infty\backslash\mathbf{Sl_2}(\mathbf{Z}), \mathbf{h'} \in \mathbf{Sl_2}(\mathbf{Z})_\infty^+\backslash\mathbf{Sl_2}(\mathbf{Z})\right\} \qquad (2.24)$$

where

$$d : \begin{cases} Gl_2(\mathbf{R}) & \longrightarrow & Sp_2(\mathbf{R}) \\ A & \longmapsto & \begin{pmatrix} (A^t)^{-1} & 0 \\ 0 & A \end{pmatrix} \end{cases}$$

$\mathbf{J} = \begin{pmatrix} 0 & -1 \\ 1 & 0 \end{pmatrix}$ and $M$ is an element of $Sl_2(\mathbf{Z})$ with $M = \begin{pmatrix} * & * \\ u & v \end{pmatrix}$. Among (2.22), (2.23),(2.24) precisely the following elements have their last row congruent to $(0, 0, *, *)$ modulo N:

$$\{ \iota_{1,1}(\mathbf{1}_2, h) \mid h \in \Gamma_\infty \backslash \Gamma_0(N) \} \tag{2.22'}$$

$$\left\{ d(\mathbf{J}) \iota_{\mathbf{1,1}}(\mathbf{h}, \mathbf{1_2}) \mid \mathbf{h} \in \Gamma_\infty^+ \backslash \Gamma_\mathbf{0}(\mathbf{N}) \right\} \tag{2.23'}$$

$$\left\{ d(M) \circ \iota_{1,1}(h, h') \mid \begin{array}{l} h \in Sl_2(\mathbf{Z})_\infty \backslash \mathbf{Sl_2(\mathbf{Z})}, \mathbf{h}' \in \mathbf{Sl_2(\mathbf{Z})_\infty^+} \backslash \mathbf{Sl_2(\mathbf{Z})} \\ M = \begin{pmatrix} * & * \\ u & v \end{pmatrix} \in Sl_2(\mathbf{Z})_\infty^+ \backslash \mathbf{Sl_2(\mathbf{Z})} \\ u, v \in \mathbf{N}, \mathbf{uc} \equiv \mathbf{0(N)}, \mathbf{vc}' \equiv \mathbf{0(N)} \end{array} \right\} \tag{2.24'}$$

Here $c$ (and $c'$) denote the lower left entry of $h$ (and $h'$).

At this point we should emphasize that (2.22')-(2.24') do **not** give representatives of

$$C_{2,1}(N) \backslash \Gamma_0^2(N),$$

since these elements are in general not in $\Gamma_0^2(N)$, but they are equivalent modulo $C_{2,1}(\mathbf{Z})$ to such representatives (by using a suitable transformation we shall finally transport them into $\Gamma_0^2(N)$).

To describe (2.24') more appropriately we fix two decompositions

$$N = N_1 \cdot N_2 \qquad \text{and} \qquad N = N_1' \cdot N_2'$$

and consider for the moment only those $h = \begin{pmatrix} * & * \\ c & d \end{pmatrix}$ and $h' = \begin{pmatrix} * & * \\ c' & d' \end{pmatrix}$ with $\gcd(c, N) = N_1$ and $\gcd(c', N) = N_1'$. Then the data

$$u, v, \begin{pmatrix} * & * \\ c & d \end{pmatrix}, \begin{pmatrix} * & * \\ c' & d' \end{pmatrix}$$

describe an element of (2.24') iff $N_2 | u$ and $N_2' | v$. These data exist only if $N | N_1 \cdot N_1'$, because we require $u, v$ to be coprime. It is a standard procedure to translate these considerations into more group theoretic terms: We denote by $\tau_{N_1} = \tau_{N_1}^N$ an element of $Sl_2(\mathbf{Z})$ with $\tau_{N_1} = \begin{pmatrix} \alpha & \beta \\ N_1 & N_2 \end{pmatrix}$ and $N_2 | \alpha$. This implies in particular that $(\tau_{N_1})^2 \in \Gamma_o(N)$.

Then (2.24') can also be described by

$$\bigcup_{\substack{N_1 N_2 = N \\ N_1' N_2' = N \\ N | N_1 N_1'}} \left\{ d(M) \iota_{1,1}(h, h') \; \middle| \; \begin{array}{l} h \in \left(\tau_{N_1}\Gamma_0(N)\tau_{N_1}^{-1}\right)_\infty \backslash \tau_{N_1}\Gamma_0(N) \\[6pt] h \in \left(\tau_{N_1'}\Gamma_0(N)\tau_{N_1'}^{-1}\right)_\infty^+ \backslash \tau_{N_1'}\Gamma_0(N) \\ u, v \text{ positive, coprime }, N_2 | u, N_2' | v \end{array} \right\} \qquad (2.24'')$$

By a routine matrix calculation, we see that (with $h, h', M, N_1, N_1'$ as in (2.24"))

$$\iota_{1,1}(\tau_{(N_1, N'_1)}, \mathbf{1_2}) \circ d(M) \circ \iota_{1,1}(h, h') \qquad (2.24''')$$

is indeed in $\Gamma_0^2(N)$, if we require (as we are allowed to do !) that $M$ is of type $M = \begin{pmatrix} r & s \\ u & v \end{pmatrix}$ with $v | r$.

We denote by $I_{m,\nu_1,\nu_2}(z_2, z_3, s)$ the $X_2^{\nu_2}X_3^{\nu_3}$−component of the $\mathbf{C[X_2, X_3]_b}$-valued function $I_m(s)$, restricted to $(z_2, z_3) \in \mathbf{H} \times \mathbf{H} \hookrightarrow \mathbf{H_2}$. With $f, \phi, \psi$ as before ($f$ now again an arbitrary element in $[\Gamma_0(N), k_1]_0$) we consider the double integral

$$\int\int_{(\Gamma_0(N) \backslash \mathbf{H})^2} \overline{\phi(z_2)\psi(z_3)} \mathcal{K}_f^{\nu_2,\nu_3}(z_2, z_3, s) y_2^{k_2} y_3^{k_3} \frac{dx_2 dy_2}{y_2^2} \frac{dx_3 dy_3}{y_3^2} \qquad (2.25)$$

where

$$\mathcal{K}_f(Z, s) =$$
$$\sum_{g \in C_{2,1}(N) \backslash \Gamma_0^2(N)} f(g < Z >^*) \left(\sigma_b \otimes \det^{k_1+a}(j(g, Z))\right)^{-1} X_2^b \cdot \left(\frac{\det \Im(g < Z >)}{\Im(g < Z >^*)}\right)^s$$

is the same Klingen-type Eisenstein series as in (2.11), but with the Hecke operator removed; again $\mathcal{K}_f^{\nu_2,\nu_3}$ denotes the $X_2^{\nu_2}X_3^{\nu_3}$-component of $\mathcal{K}_f$, restricted to $\mathbf{H} \times \mathbf{H}$. We split $\mathcal{K}_f^{\nu_2,\nu_3}$ into three parts according to the three types (2.22),(2.23),(2.24) of left cosets:

$$\mathcal{K}_f^{\nu_2,\nu_3} = \sum_{i=1}^3 \mathcal{K}_{f,i}^{\nu_2,\nu_3} \qquad (2.26)$$

It is again easy to see that $\mathcal{K}_{f,1}^{\nu_2,\nu_3}$ and $\mathcal{K}_{f,2}^{\nu_2,\nu_3}$ do not contribute to the integral (2.25). Using (2.24") we further split $\mathcal{K}_{f,3}^{\nu_2,\nu_3}$ as

$$\mathcal{K}_{f,3}^{\nu_2,\nu_3} = \sum_{N_1, N_1'} \mathcal{K}_{f, N_1, N_1'}^{\nu_2,\nu_3} \qquad (2.27)$$

We can express $\mathcal{K}_{f,N_1,N_1'}^{\nu_2,\nu_3}$ more explicitely as follows

$$\mathcal{K}_{f,N_1,N_1'}^{\nu_2,\nu_3} = 2\binom{b}{\nu_2} \sum_{h,h',u,v} f \mid_{k_1} \tau_{(N_1,N_1')}(v^2 h(<z_2>) + u^2 h'(<z_3>))$$

$$\times v^{\nu_2} u^{\nu_3} j(h,z_2)^{-k_1-a-\nu_2} j(h',z_3)^{-k_1-a-\nu_3} \times \left( \frac{\Im(h<z_2>) \cdot \Im(h'<z_3>)}{v^2 \cdot \Im(h<z_2>) + u^2 \cdot \Im(h'<z_3>)} \right)^s \quad (2.28)$$

where the summation over $h$, $h'$, $u$, $v$ is given by (2.24") with $N_1$ and $N_1'$ fixed.

It is well known how to unfold integrals like

$$\mathcal{I}(f,\phi,\psi,N_1,N_1',s) := \int_{(\Gamma_0(N)\backslash\mathbf{H})^2} \overline{\phi(z_2)\psi(z_3)} \mathcal{K}_{f,N_1,N_1'}^{\nu_2,\nu_3}(z_2,z_3,s) y_2^{k_2} y_3^{k_3} \frac{dx_2 dy_2}{y_2^2} \frac{dx_3 dy_3}{y_3^2} \quad (2.29)$$

by applying $\tau_{N_1}^{-1}$ and $\tau_{N_1'}^{-1}$: The result is

$$\mathcal{I}(f,\phi,\psi,N_1,N_1',s)$$
$$= 2\binom{b}{\nu_2} \int_{\left(\tau_{N_1}\Gamma_0(N)\tau_{N_1}^{-1}\right)_\infty \backslash \Gamma_0(N)} \int_{\left(\tau_{N_1'}\Gamma_0(N)\tau_{N_1'}^{-1}\right)_\infty \backslash \Gamma_0(N)} \overline{\left(\phi\mid_{k_2}\tau_{N_1}^{-1}\right)(z_2) \cdot \left(\psi\mid_{k_3}\tau_{N_1'}^{-1}\right)(z_3)}$$

$$\times \sum_{u,v} f \mid_{k_1} \tau_{(N_1,N_1')}(v^2 z_2 + u^2 z_3) v^{\nu_2} u^{\nu_3} \left( \frac{y_2 y_3}{v^2 y_2 + u^2 y_3} \right)^s y_2^{k_2} y_3^{k_3} \frac{dx_2 dy_2}{y_2^2} \frac{dx_3 dy_3}{y_3^2} \quad (2.30)$$

We do not want to work with Fourier expansions at several cusps, therefore we assume from now on that $f$, $\phi$, $\psi$ are normalized newforms (eigenforms of all Hecke operators) of levels $N_f$, $N_\phi$ and $N_\psi$ (all dividing $N$).

We decompose $N_1$ and $N_2$ as

$$N_1 = N_{1,f} \cdot N_1^f, \; N_2 = N_{2,f} \cdot N_2^f \quad (2.31)$$

(and the same for $N_1'$, $N_2'$ and also for $\phi$ and $\psi$)

We mention here the following facts, which we shall use in the sequel:

- $$R|N \Longrightarrow \tau_{N_1}^N = \gamma \circ \tau_{(R,N_1)}^R \quad \text{with} \quad \gamma \in \Gamma_0(R)$$

- For any divisor $d$ of $N^f$ we have

$$\begin{pmatrix} d & 0 \\ 0 & 1 \end{pmatrix} \circ \tau_{N_1}^N = \gamma \circ \begin{pmatrix} d & 0 \\ 0 & 1 \end{pmatrix} \circ \tau_{(dN_f,N_1)}^{dN_f} = \gamma \circ \tau_{N_{1,f}}^{N_f} \circ \begin{pmatrix} (d,N_1^f) & 0 \\ 0 & \frac{d}{(d,N_1^f)} \end{pmatrix}$$

with $\gamma \in \Gamma_0(N_f)$

- $$\tau_{N_1}^N = \gamma \circ \left( \prod_{q|N_2} \mathbf{V_q^N} \right) \circ \begin{pmatrix} \frac{1}{N_2} & 0 \\ 0 & 1 \end{pmatrix}$$

with $\gamma \in \Gamma_0(N)$ and (using the same notation as in [18]

$$\mathbf{V_q^N} = \left( \begin{array}{cc} x & y \\ N & q \end{array} \right)$$

with $xq - Ny = q$ and $q|x$.

- For any newform $g$ of level $N$ and weight $k$ we have (see [18, Theorem 3])

$$g \mid_k \mathbf{V_q^N} = -\mathbf{a_g(q)q^{1-\frac{k}{2}}} \cdot \mathbf{g}$$

At this point we introduce - as we already did for $N^f$ in the previous subsection- divisors $D^\phi$ and $D^\psi$ of $N^\phi$ and $N^\psi$; as usual we further factorize them (for given $N_1$ and $N_1'$) as $D^\phi = D_1^\phi \cdot D_2^\phi$ and $D^\psi = D_1'^\psi \cdot D_2'^\psi$.

Using these facts we get the following Fourier expansions

$$\left( \phi|_{k_2} \left( \begin{array}{cc} D^\phi & 0 \\ 0 & 1 \end{array} \right) |_{k_2} \tau_{N_1}^{-1} \right) (z_2) \tag{2.32}$$

$$= \phi|_{k_2} \left( \begin{array}{cc} D^\phi & 0 \\ 0 & 1 \end{array} \right) |_{k_2} \tau_{N_1}(z_2)$$

$$= \phi|_{k_2} \tau_{N_{1,\phi}}^{N_\phi} |_{k_2} \left( \begin{array}{cc} D_1^\phi & 0 \\ 0 & D_2^\phi \end{array} \right) (z_2)$$

$$= \phi|_{k_2} \left( \prod_{q|N_{2,\phi}} \mathbf{V_q^{N_\phi}} \right) |_{k_2} \left( \begin{array}{cc} \frac{D_1^\phi}{N_{2,\phi}} & 0 \\ 0 & D_2^\phi \end{array} \right) (z_2)$$

$$= \left( \prod_{q|N_{2,\phi}} -a(q)q^{-\frac{k_2}{2}+1} \right) \cdot N_{2,\phi}^{-\frac{k_2}{2}} D_1^{\phi\frac{k}{2}} D_2^{\phi\frac{-k}{2}} \cdot \phi\left( \frac{D_1^\phi}{D_2^\phi N_{2,\phi}} \cdot z_2 \right)$$

$$= D_1^{\phi\frac{k_2}{2}} D_2^{\phi-\frac{k_2}{2}} \left( \prod_{q|N_{2,\phi}} -a_\phi(q)q^{-k_2+1} \right) \sum a_\phi(n') e^{2\pi i \frac{D_1^\phi}{D_2^\phi N_{2,\phi}} \cdot n' z_2}$$

$$\left( \psi|_{k_3} \left( \begin{array}{cc} D^\psi & 0 \\ 0 & 1 \end{array} \right) |_{k_3} \tau_{N_1'}^{-1} \right) (z_3)$$

$$= D_1'^{\psi\frac{k_3}{2}} D_2'^{\psi-\frac{k_3}{2}} \left( \prod_{q|N_{2,\psi}'} -a_\psi(q)q^{-k_3+1} \right) \sum a_\psi(n'') e^{2\pi i \frac{D_1'^\psi}{D_2'^\psi N_{2,\psi}'} n'' z_3} \tag{2.33}$$

and for a divisor $d$ of $N^f$:

$$
\begin{aligned}
f|_{k_1} &\left( \begin{array}{cc} d & 0 \\ 0 & 1 \end{array} \right) \circ \tau^N_{(N_1, N_1')}(z) \\
&= f|_{k_1} \tau^{N_f}_{(N_1, N_1')_f} \circ \left( \begin{array}{cc} (d, (N_1, N_1')^f) & 0 \\ 0 & \frac{d}{(d, (N_1, N_1')^f)} \end{array} \right) \\
&= f|_{k_1} \tau^{N_f}_{(N_{1,f}, N_{1,f}')} \circ \left( \begin{array}{cc} (d, N_1^f, N_1'^f) & 0 \\ 0 & \frac{d}{(d, N_1^f, N_1'^f)} \end{array} \right) \\
&= d^{-\frac{k_1}{2}} \cdot (d, N_1^f, N_1'^f)^{k_1} \prod_{q | \operatorname{lcm}(N_{2,f}, N_{2,f}')} \left( -a_f(q) q^{-k_1+1} \right) \sum a_f(n) e^{2\pi i B \cdot n z}
\end{aligned}
\tag{2.34}
$$

We use here the simple fact that $\frac{N^f}{(N_{1,f}, N_{1,f}')} = \operatorname{lcm}(N_{2,f}, N_{2,f}')$ and

$$
B := \frac{1}{\operatorname{lcm}(N_{2,f}, N_{2,f}')} \cdot \frac{(d, N_1^f, N_1'^f)^2}{d}
\tag{2.35}
$$

Now we are ready to plug these Fourier expansions into the expression (2.30) for

$$
\mathcal{I}(f|_{k_1} \left( \begin{array}{cc} d & 0 \\ 0 & 1 \end{array} \right), \phi|_{k_2} \left( \begin{array}{cc} D^\phi & 0 \\ 0 & 1 \end{array} \right), \psi|_{k_3} \left( \begin{array}{cc} D^\psi & 0 \\ 0 & 1 \end{array} \right), N_1, N_1', s)
$$

By integration over $x_2 \bmod N_2$ and $x_3 \bmod N_2'$ we see that only those terms $a_f(n) a_\phi(n') a_\psi(n'')$ give non-zero contributions, for which

$$
B n v^2 = \frac{n'}{N_{2,\phi}} \cdot \frac{D_1^\phi}{D_2^\phi}
\tag{2.36}
$$

and

$$
B n u^2 = \frac{n''}{N_{2,\psi}'} \cdot \frac{D_1'^\psi}{D_2'^\psi}
\tag{2.37}
$$

Using

$$
\int_0^\infty \int_0^\infty e^{-4\pi B n v^2 y_2 - 4\pi B n u^2 y_3} \left( \frac{y_2 y_3}{v^2 y_2 + u^2 y_3} \right)^s y_2^{k_2-2} y_3^{k_3-2} dy_2 dy_3
$$

$$
= \frac{\Gamma(s+k_2+k_3-2)\Gamma(s+k_2-1)\Gamma(s+k_3-1)}{\Gamma(2s+k_2+k_3-2)} \cdot (4\pi B n)^{-s-k_2-k_3+2} \cdot v^{-2s-2k_2+2} \cdot u^{-2s-2k_3+2}
\tag{2.38}
$$

we obtain

$$\mathcal{I}(f|_{k_1}\begin{pmatrix} d & 0 \\ 0 & 1 \end{pmatrix}, \phi|_{k_2}\begin{pmatrix} D^\phi & 0 \\ 0 & 1 \end{pmatrix}, \psi|_{k_3}\begin{pmatrix} D^\psi & 0 \\ 0 & 1 \end{pmatrix}, N_1, N_1', s)$$

$$= 2\binom{b}{\nu_2} \cdot \mathcal{I}_\infty(s) N_2 N_2' B^{-s-k_2-k_3+2} \cdot D_1^{\phi \frac{k_2}{2}} D_2^{\phi-\frac{k_2}{2}} D_1'^{\psi \frac{k_3}{2}} D_2'^{\psi-\frac{k_3}{2}} \cdot d^{-\frac{k_1}{2}} \cdot (d, N_1^f, N_1'^f)^{k_1}$$

$$\times \prod_{q|\mathrm{lcm}(N_{2,f}, N_{2,f}')} (-a_f(q) q^{-k_1+1}) \prod_{q|N_{2,\phi}} (-a_\phi(q) q^{-k_2+1}) \prod_{q|N_{2,\psi}'} (-a_\psi(q) q^{-k_3+1}) \tag{2.39}$$

$$\times \sum_{u,v,n} a_f(n) a_\phi(nv^2 BN_{2,\phi} \cdot \frac{D_2^\phi}{D_1^\phi}) a_\psi(nu^2 BN_{2,\psi}' \cdot \frac{D_2'^\psi}{D_1'^\psi})$$

$$\times v^{-2s-2k_2+2+\nu_2} u^{-2s-2k_3+2+\nu_3} n^{-s-k_2-k_3+2}$$

with

$$\mathcal{I}_\infty(s) = (4\pi)^{-s-k_2-k_3+2} \cdot \frac{\Gamma(s+k_2+k_3-2)\Gamma(s+k_2-1)\Gamma(s+k_3-1)}{\Gamma(2s+k_2+k_3-2)} \tag{2.40}$$

To finish this section, we collect all the information obtained so far; we must take care of the fact that we worked in section (2.2) with the Eisenstein series $\mathcal{K}_f(Z, s)$ rather than with $\mathcal{K}_f(Z, s-a')$ as is required by (2.11). The value of the threefold integral $\mathcal{A}$ as defined by (2.1) is

$$\mathcal{A}(f|_{k_1}\begin{pmatrix} D^f & 0 \\ 0 & 1 \end{pmatrix}, \phi|_{k_2}\begin{pmatrix} D^\phi & 0 \\ 0 & 1 \end{pmatrix}, \psi|_{k_3}\begin{pmatrix} D^\psi & 0 \\ 0 & 1 \end{pmatrix}, s)$$

$$= 2\binom{b}{\nu_2} \cdot A(r+s, b) \cdot \mu(r+a+b, s-a') \cdot \mathcal{I}_\infty(s-a') \cdot N^{-2s-r}$$

$$\times \prod_{p\nmid N} \frac{(1-p^{-2s-r})(1-p^{2-4s-2r})}{1-p^{-2s-r+1}} \prod_{p\nmid N_f} \frac{1}{(1-\alpha_p^2 p^{-2s-r-k_1+2})(1-\alpha_p'^2 p^{-2s-r-k_1+2})} \times \prod_{p|N_f} \frac{-1}{1-p^{-2s-r}}$$

$$\times \sum_{d|N^f} \alpha(d, D^f, 2s-r) \cdot d^{-\frac{k_1}{2}} \cdot (d, N_1^f, N_1'^f)^{k_1} (D_1^\phi)^{\frac{k_2}{2}} \cdot (D_2^\phi)^{-\frac{k_2}{2}} \cdot (D_1'^\psi)^{\frac{k_3}{2}} \cdot (D_2'^\psi)^{-\frac{k_3}{2}}$$

$$\times \sum_{N_1 N_2=N, N_1' N_2'=N, N|N_1 N_1'} N_2 \cdot N_2' \cdot B^{-s+a'-k_2-k_3+2} \tag{2.41}$$

$$\times \left( \prod_{p|\mathrm{lcm}(N_{2,f}, N_{2,f}')} (-a_f(p) p^{-k_1+1}) \right) \left( \prod_{p|N_{2,\phi}} (-a_\phi(p) p^{-k_2+1}) \right) \left( \prod_{p|N_{2,\psi}'} (-a_\psi(p) p^{-k_3+1}) \right)$$

$$\times \sum_{n,u,v} a_f(n) a_\phi(nv^2 BN_{2,\phi} \cdot \frac{D_2^\phi}{D_1^\phi}) a_\psi(nu^2 BN_{2,\psi}' \cdot \frac{D_2'^\psi}{D_1'^\psi})$$

$$\times u^{-2s+a-k_3+2+\nu_3} v^{-2s+a-2k_2+2+\nu_2} n^{-s+a'-k_2-k_3+2}$$

We remind the reader that the summation over $u$ and $v$ is subject to the condition $N_2|u$ and $N_2'|v$.

### 3. The Euler factors

Using the multiplicativity properties of the Fourier coefficients of $f$, $\phi$, $\psi$ and that the conditions of summation are of multiplicative type we may now write (2.41) as

$$2 \binom{b}{\nu_2} \cdot A(r+s,b) \cdot \mu(r+a+b,s-a') \cdot \mathcal{I}_\infty(s-a') \cdot N^{-2s-r} \times \prod_p \mathcal{T}_p(s) \qquad (3.1)$$

To save notation we write

$$\mathcal{T}_p = \mathcal{T}_p^0 \cdot \mathcal{T}_p^1 \qquad (3.2)$$

where $\mathcal{T}_p^1$ denotes the last five lines of (2.41). In all cases to be considered, the pair $[(N_1)_p, N_{1p}']$ can take the values $[p,p]$, $[p,1]$ and $[1,p]$; therefore we split $\mathcal{T}_p^1$ as

$$\mathcal{T}_p^1 = C_{pp} + C_{p1} + C_{1p} \qquad (3.3)$$

If $d_p$ can take both values $1$ and $p$, we further decompose $C_{**}$ as

$$C_{**} = C_{**}^1 + C_{**}^p \qquad (3.4)$$

according to the cases $d_p = 1$ and $d_p = p$.

**Part I:** $D_p^f = D_p^\phi = D_p^\psi = 1$

**IA : The case** $p|N_f,\ p|N_\phi,\ p|N_\psi$

(This is the case also considered by Gross/Kudla[11])

The conditions imply $d_p = N_p^f = N_p^\phi = N_p^\psi = 1$, $\alpha_p(d, D^f, s) = 1$ and

$$B_p = \begin{cases} 1 & \text{for the case} \quad [p,p] \\ \frac{1}{p} & \text{otherwise} \end{cases}$$

$$C_{pp}(s) = \left( \sum_{l=0}^\infty \sum_{t=0}^\infty a_f(p^l) a_\phi(p^{l+2t}) a_\psi(p^l) \left( p^l \right)^{-s+a'-k_2-k_3+2} \left( p^t \right)^{-2s+a-2k_2+2+\nu_2} \right.$$

$$\left. + \sum_{l=0}^\infty \sum_{t=1}^\infty a_f(p^l) a_\phi(p^l) a_\psi(p^{l+2t}) \left( p^l \right)^{-s+a'-k_2-k_3+2} \left( p^t \right)^{-2s+a-2k_3+2+\nu_3} \right) \qquad (3.5)$$

This equals (we use the fact that $a_\phi(p^2)$ and $a_\psi(p^2)$ are powers of $p$)

$$\frac{1}{1-a_f(p)a_\phi(p)a_\psi(p)p^{-s+a'-k_2-k_3+2}}\cdot\left\{\frac{1}{1-a_\phi(p^2)p^{-2s+a-2k_2+2+\nu_2}}+\frac{a_\psi(p^2)p^{-2s+a-2k_3+2+\nu_3}}{1-a_\psi(p)p^{-2s+a-k_3+2+\nu_3}}\right\}$$

$$=\frac{1}{1-a_f(p)a_\phi(p)a_\psi(p)p^{-s+a'-k_2-k_3+2}}\cdot\left\{\frac{1}{p^{-2s+a-k_2+\nu_2}}+\frac{p^{-2s+a-k_3+\nu_3}}{1-p^{-2s+a-k_3+\nu_3}}\right\}$$

$$=\frac{1}{1-a_f(p)a_\phi(p)a_\psi(p)p^{-s-3a'-2r-b+2}}\cdot\frac{1+p^{-2s-r}}{1-p^{-2s-r}}.$$

$$C_{1,p}=p\cdot p^{s-a'+k_2+k_3-2}\cdot(-a_f(p)p^{-k_1+1})\cdot(-a_\phi p^{-k_2+1})\times \qquad (3.6)$$

$$\times\sum_{l=0}^{\infty}\sum_{t=1}^{\infty}a_f(p^l)a_\phi(p^l\cdot\frac{1}{p}\cdot p)a_\psi(p^{l+2t}\cdot\frac{1}{p})\left(p^t\right)^{-2s+a-2k_3+2+\nu_3}\cdot\left(p^l\right)^{-s+a'-k_2-k_3+2}$$

$$=\frac{1}{1-a_f(p)a_\phi(p)a_\psi(p)p^{-s+a'-k_2-k_3+2}}\cdot\frac{a_f(p)a_\phi(p)a_\psi(p)p^{3-s-3a'-b-2r}}{1-a_\psi(p^2)p^{-2s+a-2k_3+2+\nu_3}}$$

$$=\frac{1}{1-a_f(p)a_\phi(p)a_\psi(p)p^{-s-3a'-2r-b+2}}\cdot\frac{a_f(p)a_\phi(p)a_\psi(p)p^{3-s-3a'-b-2r}}{1-p^{-2s-r}}$$

Quite the same computation shows that

$$C_{p,1}=C_{1,p} \qquad (3.7)$$

Hence

$$\mathcal{T}_p^1=\frac{1}{1-a_f(p)a_\phi(p)a_\psi(p)p^{-s-3a'-2r-b+2}}$$
$$\times\frac{1+2a_f(p)a_\phi(p)a_\psi(p)p^{-s-3a'-b-2r+3}+p^{-2s-r}}{1-p^{-2s-r}} \qquad (3.8)$$

Now we write the numerator as

$$1+2a_f(p)a_\phi(p)a_\psi(p)p^{-s-3a'-b-2r+3}+p^{-2s-r}=\left(1+a_f(p)a_\phi(p)a_\psi(p)p^{-s+3-3a'-b-2r}\right)^2$$
$$=\frac{(1-p^{-2s-r})^2}{(1-a_f(p)a_\phi(p)a_\psi(p)p^{-s+3-3a'-b-2r})^2} \qquad (3.9)$$

Therefore we obtain

$$\mathcal{T}_p(s)=\frac{-1}{1-a_f(p)a_\phi(p)a_\psi(p)p^{-s-3a'-2r-b+2}}\cdot\frac{1}{(1-a_f(p)a_\phi(p)a_\psi(p)p^{-s+3-3a'-b-2r})^2} \qquad (3.10)$$

**IB :The case $p|N^f$, $p|N_\phi$, $p|N_\psi$** These conditions imply

$$(N_f)_p=N_p^\phi=N_p^\psi=D_p^\phi=D_p^\psi=1 \qquad\text{and}\qquad d_p\in\{1,p\}$$

$$\alpha_p(d, D^f, s) = \begin{cases} -(1 + p^{1-s}) & \text{if} \quad d_p = 1 \\ p^{1-\frac{k_1}{2}} \cdot a_f(p) & \text{if} \quad d_p = p \end{cases}$$

and

$$B_p = \begin{cases} p & \text{if} \quad d_p = (N_1)_p = N'_{1p} = p \\ \frac{1}{d_p} & \text{otherwise} \end{cases}$$

What we get is this:

$$C_{pp}^1(s) = \frac{-(1 + p^{1-2s-r})(1 + p^{-2s-r})}{(1 - \alpha_p a_\phi(p) a_\psi(p) p^{-s-3a'-2r-b+2})(1 - \alpha_p a_\phi(p) a_\psi(p) p^{-s-3a'-2r-b+2})(1 - p^{-2s-r})} \tag{3.11}$$

(similar calculation as in (3.5)

$$
\begin{aligned}
C_{1p}^1(s) &= -(1 + p^{1-2s-r}) \cdot p \cdot (-a_\phi(p) p^{-k_2+1}) \times \\
&\qquad \sum_{l=0}^{\infty} \sum_{t=1}^{\infty} a_f(p^l) a_\phi(p^l \cdot p) a_\psi(p^{l+2t})(p^t)^{-2s+a-2k_3+2+\nu_3} \cdot (p^l)^{-s+a'-k_2-k_3+2} \\
&= \frac{(1 + p^{1-2s-r}) \cdot p^{-2s-r}}{(1 - \alpha_p a_\phi(p) a_\psi(p) p^{-s-3a'-2r-b+2})(1 - \alpha_p a_\phi(p) a_\psi(p) p^{-s-3a'-2r-b+2})(1 - p^{-2s-r})}
\end{aligned}
\tag{3.12}
$$

$$C_{p1}^1 = C_{1p}^1 \tag{3.13}$$

$$
\begin{aligned}
C_{pp}^p(s) &= p^{1-\frac{k_1}{2}} a_f(p) p^{-\frac{k_1}{2}} p^{k_1} p^{-s+a'-k_2-k_3+2} \times \\
&\quad \left\{ \sum_{l=0}^{\infty} \sum_{t=0}^{\infty} a_f(p^l) a_\phi(p^{l+2t} p) a_\psi(p^l p)(p^t)^{-2s+a-2k_2+2+\nu_2}(p^l)^{-s+a'-k_2-k_3+2} + \right. \\
&\qquad \left. \sum_{l=0}^{\infty} \sum_{t=1}^{\infty} a_f(p^l) a_\phi(p^l p) a_\psi(p^{l+2t} p)(p^t)^{-2s+a-2k_3+2+\nu_3}(p^l)^{-s+a'-k_2-k_3+2} \right\} \\
&= \frac{a_f(p) a_\phi(p) a_\psi(p) p^{-s-3a'-2r-b+3}(1 + p^{-2s-r})}{(1 - \alpha_p a_\phi(p) a_\psi(p) p^{-s-3a'-2r-b+2})(1 - \alpha'_p a_\phi(p) a_\psi(p) p^{-s-3a'-2r-b+2})(1 - p^{-2s-r})}
\end{aligned}
\tag{3.14}
$$

$$
\begin{aligned}
C_{1,p}^p &= p^{1-\frac{k_1}{2}} a_f(p) p^{-\frac{k_1}{2}} \cdot p \cdot p^{s-a'+k_1+k_2-2}(-a_\phi(p) p^{-k_2+1}) \times \\
&\quad \sum_{l=0}^{\infty} \sum_{t=1}^{\infty} a_f(p^l) a_\phi(p^l \tfrac{1}{p} \cdot p) a_\psi(p^{l+2t} \tfrac{1}{p})(p^t)^{-2s+a-k_3+2+\nu_3}(p^l)^{-s+a'-k_2-k_3+2} \\
&= \frac{-a_f(p) a_\phi(p) a_\psi(p) p^{-s-2r-3a'-b+3}}{(1 - \alpha_p a_\phi(p) a_\psi(p) p^{-s-3a'-2r-b+2})(1 - \alpha'_p a_\phi(p) a_\psi(p) p^{-s-3a'-2r-b+2})(1 - p^{-2s-r})}
\end{aligned}
\tag{3.15}
$$

$$C_{p1}^p = C_{1p}^p \tag{3.16}$$

All the $C_{**}^*$ have the same denominator, the numerator of $\sum C_{**}^*$ being equal to

$$1 - p^{-2s-r} + p^{1-2s-r} - a_f a_\phi a_\psi \left\{ p^{-3s-3r-a-b+3} - p^{-s-a-2r-b+3} \right\}$$

This can be factorized as

$$-(1 - p^{-2s-r})(1 + \alpha_p a_\phi(p) a_\psi(p) p^{-s-2r-3a'-b+3})(1 + \alpha'_p a_\phi(p) a_\psi(p) p^{-s-2r-3a'-b+3})$$

Now we also write the denominator of $\mathcal{T}_p^0$ as product of linear factors in $p^{-s}$ using

$$(1 - \alpha_p^2 p^{-k_1+2-2s-r}) = (1 - \alpha_p a_\phi a_\psi p^{-s-3a'-b-2r+3})(1 + \alpha_p a_\phi a_\psi p^{-s-3a'-b-2r+3})$$

(and the same with $\alpha'_p$ instead of $\alpha_p$); therefore we get some cancellations and arrive at

$$
\begin{aligned}
\mathcal{T}_p &= \frac{-1}{(1 - \alpha_p a_\phi(p) a_\psi(p) p^{-s-3a'-b-2r+2})(1 - \alpha'_p a_\phi(p) a_\psi(p) p^{-s-3a'-b-2r+2})} \\
&\times \frac{1}{(1 - \alpha_p a_\phi(p) a_\psi(p) p^{-s-3a'-b-2r+3})(1 - \alpha'_p a_\phi(p) a_\psi(p) p^{-s-3a'-b-2r+3})}
\end{aligned}
\tag{3.17}
$$

**IC The case** $p|N_f$, $p|N^\phi$, $p|N^\psi$

We have

$$(B)_p = \begin{cases} 1 & \text{if} \quad N_1 = N'_1 = p \\ \frac{1}{p} & \text{otherwise} \end{cases}$$

and we get

$$
\begin{aligned}
C_{1p} &= p \cdot p^{s-a'+k_2+k_3-2} \cdot (-a_f(p) p^{-k_1+1}) \times \\
&\quad \sum_{l=1}^\infty \sum_{t=1}^\infty a_f(p^l) a_\phi(p^l \tfrac{1}{p}) a_\psi(p^{l+2t} \tfrac{1}{p})(p^t)^{-2s+a-2k_3+2+\nu_3}(p^l)^{-s+a'-k_2-k_3+2} \\
&= -\sum_{l'=0}^\infty \sum_{t=1}^\infty a_f(p^{l'}) a_\phi(p^{l'}) a_\psi(p^{l'+2t})(p^t)^{-2s+a-2k_3+2+\nu_3}(p^{l'})^{-s+a'-k_2-k_3+2}
\end{aligned}
\tag{3.18}
$$

$$
C_{p1} = -\sum_{l=0}^\infty \sum_{t=1}^\infty a_f(p^l) a_\phi(p^{l+2t}) a_\psi(p^l)(p^t)^{-2s+a-2k_2+2+\nu_2}(p^l)^{-s+a'-k_2-k_3+2}
\tag{3.19}
$$

$$
C_{pp}(s) = \left( \sum_{l=0}^\infty \sum_{t=0}^\infty a_f(p^l) a_\phi(p^{l+2t}) a_\psi(p^l) \left( p^l \right)^{-s+a'-k_2-k_3+2} \left( p^t \right)^{-2s+a-2k_2+2+\nu_2} + \right.
\tag{3.20}
$$

$$
\left. + \sum_{l=0}^\infty \sum_{t=1}^\infty a_f(p^l) a_\phi(p^l) a_\psi(p^{l+2t}) \left( p^l \right)^{-s+a'-k_2-k_3+2} \left( p^t \right)^{-2s+a-2k_3+2+\nu_3} \right)
$$

Summing up the $C_{**}$ we see that the summands in (3.20) with $t \neq 0$ cancel against (3.18) and (3.19) and we get

$$
\begin{aligned}
\mathcal{T}_p(s) &= -(1 - p^{-2s-r})^{-1} \sum_{l=0}^{\infty} a_f(p^l) a_\phi(p^l) a_\psi(p^l) (p^l)^{-s-3a'-b-2r+2} \\
&= \frac{-1}{(1 - a_f(p)\beta_p\gamma_p p^{-s-3a'-b-2r+2})(1 - a_f(p)\beta_p\gamma_p' p^{-s-3a'-b-2r+2})} \times \\
&\quad \times \frac{1}{(1 - a_f(p)\beta_p'\gamma_p p^{-s-3a'-b-2r+2})(1 - a_f(p)\beta_p'\gamma_p' p^{-s-3a'-b-2r+2})}
\end{aligned}
\tag{3.21}
$$

### I D: The case $p \nmid N$

This case (which is in some sense the most difficult one) was previously considered in [9] for the case of equal weights. It was already noticed in [21] and [22] that the result from [9] carries over to the case of arbitrary weights. We just state the result here:

$$
\mathcal{T}_p = (1 - p^{-2s-r})(1 - p^{2-4s-2r}) L_p(f \otimes \phi \otimes \psi, s + 3a' + 2r + b - 2)
\tag{3.22}
$$

### Part II: $D_p^f = N_p^f$ , $D_p^\phi = N_p^\phi$, $D_p^\psi = N_p^\psi$

Only the cases B and C are to be investigated, the case B being the most complicated:

### II B : The case $p|N^f$, $p|N_\phi$, $p|N_\psi$

These conditions imply

$$
\alpha_p(d, D^f, s) = \begin{cases} -a_f(p)p^{2-\frac{k_1}{2}-s} & \text{if} \quad d_p = 1 \\ p(1 + p^{1-s}) & \text{if} \quad d_p = p \end{cases}
$$

$$
B_p = \begin{cases} p & \text{if} \quad d_p = (N_1)_p = (N_1')_p = p \\ \frac{1}{d} & \text{otherwise} \end{cases}
$$

$$
\begin{aligned}
C_{pp}^1 &= \alpha_p(1, D^f, 2s+r) \times \\
&\quad \Big\{ \sum_{l=0}^{\infty} \sum_{t=0}^{\infty} a_f(p^l) a_\phi(p^{l+2t}) a_\psi(p^l)(p^l)^{-2s+a-2k_2+2+\nu_2} (p^t)^{-s+a'-k_2-k_3+2} + \\
&\quad \sum_{l=0}^{\infty} \sum_{t=1}^{\infty} a_f(p^l) a_\phi(p^l) a_\psi(p^{l+2t})(p^t)^{-2s+a-2k_3+2+\nu_3} (p^l)^{-s+a'-k_2-k_3+2} \Big\} \\
&= \frac{\alpha_p(1, D^f, s)(1 + p^{-2s-r})}{(1 - \alpha_p a_\phi(p) a_\psi(p) p^{-s-3a'-2r-b+2})(1 - \alpha_p' a_\phi(p) a_\psi(p) p^{-s-3a'-2r-b+2})(1 - p^{-2s-r})}
\end{aligned}
\tag{3.23}
$$

(computation similar to (3.5))

$$
\begin{aligned}
C_{1p}^1 &= \alpha_p(1, D^f, s) \cdot p \cdot (-a_\phi(p)p^{-k_2+1}) \times \\
&\quad \sum_{l=0}^{\infty}\sum_{t=1}^{\infty} a_f(p^l)a_\phi(p^l p)a_\psi(p^{l+2t})(p^t)^{-2s+a-2k_3+2+\nu_3}(p^l)^{-s+a'-k_2-k_3+2} \\
&= \frac{-p^{-2s-r}\alpha_p(1, D^f, 2s+r)}{(1-\alpha_p a_\phi(p)a_\psi(p)p^{-s-3a'-2r-b+2})(1-\alpha_p' a_\phi(p)a_\psi(p)p^{-s-3a'-2r-b+2})(1-p^{-2s-r})}
\end{aligned}
\tag{3.24}
$$

$$
C_{p,1}^1 = C_{1,p}^1 \tag{3.25}
$$

and

$$
\sum C_{**}^1 = \frac{-a_f(p)p^{2-\frac{k_1}{2}-s}}{(1-\alpha_p a_\phi(p)a_\psi(p)p^{-s-3a'-2r-b+2})(1-\alpha_p' a_\phi(p)a_\psi(p)p^{-s-3a'-2r-b+2})} \tag{3.26}
$$

$$
\begin{aligned}
C_{pp}^p(s) &= \alpha_p(p, D^f, s+2r) \cdot p^{-\frac{k_1}{2}}p^{k_1}p^{-s+a'-k_2-k_3+2} \times \\
&\quad \sum_{l=0}^{\infty}\sum_{t=0}^{\infty} a_f(p^l)a_\phi(p^{l+2t}p)a_\psi(p^l p)(p^t)^{-2s+a-2k_2+2+\nu_2}(p^l)^{-s+a'-k_2-k_3+2} + \\
&\quad \sum_{l=0}^{\infty}\sum_{t=1}^{\infty} a_f(p^l)a_\phi(p^l p)a_\psi(p^{l+2t}p)(p^t)^{-2s+a-2k_3+2+\nu_3}(p^l)^{-s+a'-k_2-k_3+2} \\
&= \frac{\alpha_p(p, D^f, 2s+r)a_\phi(p)a_\psi(p)p^{\frac{k_1}{2}-s-3a'-2r-b+2}(1+p^{-2s-r})}{(1-\alpha_p a_\phi(p)a_\psi(p)p^{-s-3a'-2r-b+2})(1-\alpha_p a_\phi(p)a_\psi(p)p^{-s-3a'-2r-b+2})(1-p^{-2s-r})}
\end{aligned}
\tag{3.27}
$$

(similar computation as in (3.5))

$$
\begin{aligned}
C_{1p}^p(s) &= \alpha_p(p, D^f, 2s+r) \cdot p^{-\frac{k_1}{2}} \cdot p \cdot p^{s-a'+k_2+k_3-2}(-a_\phi p^{-k_2+1}) \times \\
&\quad \sum_{l=0}^{\infty}\sum_{t=0}^{\infty} a_f(p^l)a_\phi(p^l\tfrac{1}{p}p)a - \psi(p^{l+2t}\tfrac{1}{p})(p^t)^{-2s+a-2k_3+2+\nu_3}(p^l)^{-s+a'-k_2-k_3+2} \\
&= \frac{-\alpha_p(p, D^f, 2s+r)a_\phi(p)a_\psi(p)p^{-\frac{k_1}{2}-s-a'-r+2}}{(1-\alpha_p a_\phi(p)a_\psi(p)p^{-s-3a'-2r-b+2})(1-\alpha_p' a_\phi(p)a_\psi(p)p^{-s-3a'-2r-b+2})(1-p^{-2s-r})}
\end{aligned}
\tag{3.28}
$$

$$
C_{p1}^p = C_{1p}^p \tag{3.29}
$$

and hence

$$
\sum C_{**}^p = \frac{-a_\phi a_\psi(1+p^{1-2s-r})p^{-\frac{k_1}{2}-s-a'-r+3}}{(1-\alpha_p a_\phi(p)a_\psi(p)p^{-s-3a'-2r-b+2})(1-\alpha_p' a_\phi(p)a_\psi(p)p^{-s-3a'-2r-b+2})} \tag{3.30}
$$

The numerator of $\sum C_*^*$ is equal to

$$-\frac{1}{a_\phi a_\psi}p^{-\frac{k_1}{2}+r+3a'+b-1-s}(1+a_f(p)a_\phi(p)a_\psi p^{-s-2r+3-b-3a'}+p^{1-2s-r}) \tag{3.31}$$

We can therefore apply the same kind of trick as in case I B; denoting by $\epsilon_p(\phi)$ the eigenvalue of the Atkin-Lehner-involution $V_p^N$ acting on $\phi$ (and similarly for $\psi$) we end up with

$$
\begin{aligned}
\mathcal{T}_p(s) &= \frac{-\epsilon_p(\phi)\epsilon_p(\psi)p^{1-\frac{r}{2}-s}}{(1-\alpha_p a_\phi(p)a_\psi(p)p^{-s-3a'-b-2r+2})(1-\alpha'_p a_\phi(p)a_\psi(p)p^{-s-3a'-b-2r+2})} \\
&\times \frac{1}{(1-\alpha_p a_\phi(p)a_\psi(p)p^{-s-3a'-b-2r+3})(1-\alpha'_p a_\phi(p)a_\psi(p)p^{-s-3a'-b-2r+3})}
\end{aligned} \tag{3.32}
$$

## II C: The case $p|N_f,\ p|N^\phi,\ p|N^\psi$

$$
\begin{aligned}
C_{pp}(s) &= p^{\frac{k_2}{2}}\cdot p^{\frac{k_3}{2}}\times \\
&\quad \Big\{\sum_{l=1}^\infty \sum_{t=0}^\infty a_f(p^l)a_\phi(p^{l+2t}\tfrac{1}{p})a_\psi(p^l\tfrac{1}{p})(p^t)^{-2s+a-2k_2+2+\nu_2}(p^l)^{-s+a'-k_2-k_3+2}+ \\
&\quad +\sum_{l=1}^\infty \sum_{t=1}^\infty a_f(p^l)a_\phi(p^l\tfrac{1}{p})a_\psi(p^{l+2t}\tfrac{1}{p})(p^t)^{-2s+a-2k_3+2+\nu_3}(p^l)^{-s+a'-k_2-k_3+2}\Big\} \\
&= p^{-s+a'+2-\frac{k_2+k_3}{2}}\times \\
&\quad \Big\{\sum_{l=0}^\infty \sum_{t=0}^\infty a_f(p^{l+1})a_\phi(p^{l+2t})a_\psi(p^l)(p^t)^{-2s+a-2k_2+2+\nu_2}(p^l)^{-s+a'-k_2-k_3+2}+ \\
&\quad \sum_{l=0}^\infty \sum_{t=1}^\infty a_f(p^{l+1})a_\phi(p^l)a_\psi(p^{l+2t})(p^t)^{-2s+a-2k_3+2+\nu_3}(p^l)^{-s+a'-k_2-k_3+2}\Big\}
\end{aligned} \tag{3.33}
$$

$$
\begin{aligned}
C_{1p} &= p^{-\frac{k_2}{2}+\frac{k_3}{2}}\cdot p\cdot p^{s-a'+k_2+k_3-2}\cdot(-a_f(p)p^{-k_1+1})\times \\
&\quad \sum_{l=0}^\infty \sum_{t=1}^\infty a_f(p^l)a_\phi(p^l\tfrac{1}{p})a_\psi(p^{l+2t}\tfrac{1}{p}\tfrac{1}{p})(p^t)^{-2s+a-2k_3+2+\nu_3}(p^l)^{-s+a'-k_2-k_3+2} \\
&= -p^{-s+a'+2-\frac{k_2+k_3}{2}}\times \\
&\quad \sum_{l=0}^\infty \sum_{t=0}^\infty a_f(p^{l+1})a_\phi(p^l)a_\psi(p^{l+2t})(p^t)^{-2s+a-2k_3+2+\nu_3}(p^l)^{-s+a'-k_2-k_3+2}
\end{aligned} \tag{3.34}
$$

$$
\begin{aligned}
C_{p1} &= -p^{-s+a'+2-\frac{k_2+k_3}{2}}\times \\
&\quad \sum_{l=0}^\infty \sum_{t=0}^\infty a_f(p^{l+1})a_\phi(p^{l+2t})a_\psi(p^l)(p^t)^{-2s+a-2k_2+2+\nu_3}(p^l)^{-s+a'-k_2-k_2+2}
\end{aligned} \tag{3.35}
$$

Hence in the sum of the $C_{**}$ only the "$t=0$-part" of $C_{p1}$ survives and we obtain

$$
\begin{aligned}
\mathcal{T}_p &= \frac{-\epsilon_p(f)p^{1-s-\frac{r}{2}}}{(1-a_f(p)\beta_p\gamma_p p^{-s-3a'-b-2r+2})(1-a_f(p)\beta_p\gamma'_p p^{-s-3a'-b-2r+2})} \\
&\times \frac{1}{(1-a_f(p)\beta'_p\gamma_p p^{-s-3a'-b-2r+2})(1-a_f(p)\beta'_p\gamma'_p p^{-s-3a'-b-2r+2})}
\end{aligned} \tag{3.36}
$$

**Remark**: Although our list of Euler factors is complete, the reader should be aware of the fact that in our integral representation (2.1) we are free to interchange the roles of $\phi$ and $\psi$ (interchanging the roles of $\nu_2$ and $\nu_3$ at the same time), but $f$ has to be the cusp form of largest weight. Therefore e.g. in case IB, IIB we should also consider the case where $p$ does not divide the level of $\phi$ or $\psi$. It will be left to the reader to show by similar computations as above that in those cases the Euler factor will be the same (as should be more or less clear from an adelic point of view).

## 4. The Functional equation

The factors $\mathcal{T}_p = \mathcal{T}_p(s)$ computed in the previous section are for $p \nmid N$ and for $p \mid \gcd(N_f, N_\phi, N_\psi)$ up to a shift in the argument and an elementary factor the known Euler factors of the triple product $L$-function $L(f, \phi, \psi, s)$ associated to $f, \phi, \psi$. We define therefore now:

**Definition 4.1** *The triple product $L$-function associated to $f, \phi, \psi$ is for $s > \frac{k_1+k_2+k_3-1}{2}$ defined as $L(f, \phi, \psi, s) = \prod_p L_p(f, \phi, \psi, s)$ where the Euler factor $L_p(f, \phi, \psi, s)$ is given by*

$$L_p(f, \phi, \psi, s + 3a' + 2r + b - s) = \begin{cases} -\mathcal{T}_p(s) & \text{if } p \mid N \\ (1 - p^{-2s-r})^{-1}(1 - p^{2-4s-2r})^{-1}\mathcal{T}_p(s) & \text{if } p \nmid N \end{cases} \quad (4.1)$$

*With*

$$L_\infty(f, \phi, \psi, s) = \Gamma_{\mathbf{C}}(s)\Gamma_{\mathbf{C}}(s + 1 - k_1)\Gamma_{\mathbf{C}}(s + 1 - k_2)\Gamma_{\mathbf{C}}(s + 1 - k_3) \quad (4.2)$$

*(where $\Gamma_{\mathbf{C}}(s) = (2\pi)^{-s}\Gamma(s)$ as usual) the completed triple product $L$-function is*

$$\Lambda(f, \phi, \psi, s) = L_\infty(f, \phi, \psi, s)L(f, \phi, \psi, s).$$

With these notations we have for the integral $\mathcal{A}(f, \phi, \psi, s)$ from 2.1 (with $c_2(s)$ as in (1.19), $\zeta_p(s) = (1 - p^{-s})^{-1}$ for finite primes $p$ and $\zeta_\infty(s) = \pi^{-s/2}\Gamma(s/2)$):

**Theorem 4.2**

$$\frac{\mathcal{A}(f, \phi, \psi, s)}{c_2(s)}$$

$$= \binom{b}{\nu_2} i^{b+k_1} N^{-2s-2} 2^{3-2b-a} \frac{\pi^{2+3a'}}{(s+1)(2s+1)} \frac{L_\infty(f, \phi, \psi, s + \frac{k_1+k_2+k_3}{2} - 1)}{\zeta_\infty(4s+2)\zeta_\infty(2s+2)}$$

$$\times \prod_{p \nmid N} \frac{L_p(s + \frac{k_1+k_2+k_3}{2} - 1)}{\zeta_p(4s+2)\zeta_p(2s+2)} \prod_{p \mid N} (-L_p(f, \phi, \psi, s + \frac{k_1+k_2+k_3}{2} - 1)) \quad (4.3)$$

*Proof.* This follows from 2.41 and the results of Section 3.

**Theorem 4.3** *The function* $\Lambda(f, \phi, \psi, s)$ *has a meromorphic continuation to all of* $\mathbf{C}$ *and satisfies the functional equation*

$$\Lambda(f, \phi, \psi, s)$$
$$= -N^{-4(s+\frac{k_1+k_2+k_3}{2}-1)}(\gcd(N_f, N_\phi, N_\psi))^{-(s+\frac{k_1+k_2+k_3}{2}-1)}(\prod_{p|N}\epsilon_p)\Lambda(k_1+k_2+k_3-2-s) \quad (4.4)$$

*where for* $p \mid N$ *the number* $\epsilon_p$ *is defined as the product of the eigenvalues under the* $p$-*Atkin-Lehner-involution* $w_p$ *of those forms among* $f, \phi, \psi$ *whose level is divisible by* $p$.

*Proof.* We put $r = 2$ since for this choice the functional equation of the Eisenstein series is under $s \mapsto -s$. It is then not difficult to read off the functional equation of $\mathbf{E}$ from the calculations in [11]; actually things become somewhat simplified since we need the functional equation only up to oldforms. We need the following Lemma:

**Lemma 4.4** *The Eisenstein series* $\mathbf{E_{2,s}}$ *satisfies the functional equation:*

$$(s+1)(2s+1)\zeta_\infty(4s+2)\zeta_\infty(2s+2)\prod_{p|N}\zeta_p(4s+2)\zeta_p(2s+2)\mathbf{E_{2,s}(Z)}$$

$$= -N^{-9s}(s-1)(2s-1)\zeta_\infty(-4s+2)\zeta_\infty(-2s+2) \quad (4.5)$$

$$\times \prod_{p\nmid N}\zeta_p(-4s+2)\zeta_p(-2s+2)\mathbf{E_{2,-s}}\begin{pmatrix} 0_3 & -1_3 \\ N1_3 & 0_3 \end{pmatrix}(\mathbf{Z}) + \widetilde{\mathbf{E_s}}$$

*where* $\widetilde{\mathbf{E_s}}$ *is a linear combination of Eisenstein series for groups* $\Gamma_0^3(N')$ *with* $N'$ *strictly dividing* $N$ *or conjugates of these by a matrix* $\begin{pmatrix} 0_3 & -1_3 \\ M1_3 & 0_3 \end{pmatrix}$ *with* $M \mid (N/N')$

*Proof.* Let $\mathbf{F_{r,s}(Z)}$ be the Eisenstein series of degree 3 defined in the same way as $\mathbf{E_{r,s}}$ in Section 2, but with the summation running over coprime symmetric pairs $(C, D)$ with $\gcd(\det C, N) = 1$; one has

$$\mathbf{E_{r,s}|_r}\begin{pmatrix} 0 & -1 \\ N & 0 \end{pmatrix} = \mathbf{N^{-(3r/2)-3s}F_{r,s}} \quad (4.6)$$

The calculation of the local intertwining operators $M(s) = M_p(s)$ in sections 5 and 6 of [11] gives (in the notations of that article) $M_p(s)\Phi_p(s)$ for $p = \infty$ or $p \nmid N$ and allows for $p \nmid N$ to express $M_p(s)\Phi_p^0(s)$ explicitly as a linear combination of the sections $\Phi_p^0(-s), \Phi_p^3(-s), (\Phi_p)_K(-s), (\Phi_p)_{K'}(-s)$. An elementary calculation shows in this case that the coefficient at $\Phi_p^0(-s)$ is zero and that the coefficient at $\Phi_p^3(-s)$ equals

$$p^{-6s-3}\frac{\zeta_p(2s-1)\zeta_p(4s-1)}{\zeta_p(-4s+2)\zeta_p(-2s+2)}.$$

We let $Z = X + iY$ and $g = (g_\infty, 1, \dots) \in Sp_3(\mathbf{A})$ with $g_\infty = \begin{pmatrix} 1_3 & X \\ 0_3 & 1_3 \end{pmatrix}\begin{pmatrix} Y^{1/2} & 0_3 \\ 0_3 & Y^{-1/2} \end{pmatrix}$. Then it is well known that with

$$\Phi^0 = \Phi_\infty^{-2} \times \prod_{p|N}\Phi_p^0 \times \prod_{p\nmid N}(\Phi_p)_K$$

and analogously defined $\Phi^3$ one has $\mathbf{E_{2,s}(\bar{Z})} = (\det \mathbf{Y})^{-1}\mathbf{E(g,s,\Phi^0)}$ and

$$\mathbf{F_{2,s}(\bar{Z})} = (\det \mathbf{Y})^{-1}\mathbf{E(g,s,\Phi^3)}$$

(see Proposition 7.5 of [11]), and analogous formulae are true for all the global sections $\Phi$ for which the $p$-adic components $\Phi_p$ are one of the $\Phi_p^0, \Phi_p^3, (\Phi_p)_K, (\Phi_p)_{K'}$. The assertion of the lemma follows upon using the duplication formula for the gamma-function for the contribution from the infinite place and then applying the functional equation of the Riemann zeta function.

We can now finish the proof of Theorem 4.3. From Theorem 4.2, Lemma 4.4 and the fact that the integrand in $\frac{A(f,\phi,\psi,s)}{c_2(s)}$ contains by 1.18 a differential operator independent of $s$ we see (using the results of Part II in Section 3) that

$$\Lambda(s + \frac{k_1+k_2+k_3}{2} - 1) = -N^{-4s}(\gcd(N_f, N_\phi, N_\psi))^{-s}(\prod_{p|N}\epsilon_p)\Lambda(-s + \frac{k_1+k_2+k_3}{2} - 1) \tag{4.7}$$

where for $p \mid N$ the number $\epsilon_p$ is defined as in the assertion of Theorem 4.3. Putting $s' = s + \frac{k_1+k_2+k_3}{2} - 1$ we obtain the desired functional equation under $s' \mapsto k_1 + k_2 + k_3 - 2 - s'$.

## 5.Computation of the central critical value

Following the strategy of [12] we now evaluate the integral 2.2 at the point $s = 0$ using a variant of Siegel's theorem, i. e. , expressing the value at $s = 0$ of the Eisenstein series $G$ as a linear combination of theta series. The setting for this is basically the same as in [2, 5, 4]. Let $M_1, M_2, M_3$ be relatively prime square free integers such that $M_1$ has an odd number of prime divisors. By $D = D(M_1)$ we denote the quaternion algebra over $\mathbf{Q}$ ramified at $\infty$ and the primes dividing $M_1$ and by $R = R(M_1, M_2)$ an Eichler order of level $M = M_1 M_2$ in $D$, i. e. , the completion $R_p$ is a maximal order for $p \nmid M_2$ and conjugate to $\{\left(\begin{smallmatrix} a & b \\ c & d \end{smallmatrix}\right) \in M_2(\mathbf{Z_p}) \mid \mathbf{c} \equiv \mathbf{0} \bmod \mathbf{p}\}$ for $p \mid M_2$, where we identify $D \otimes \mathbf{Q_p}$ with $M_2(\mathbf{Q_p})$ for $p \nmid M_1$. By $\text{gen}(M_1, M_2, M_3)$ we denote the genus of $\mathbf{Z}$-lattices with quadratic form (quadratic lattices) of $R(M_1, M_2)$ equipped with the norm form of the quaternion algebra scaled by $M_3$. The genus theta series of degree $n$ of $\text{gen}(M_1, M_2, M_3)$ is then

$$\Theta_{M_1,M_2,M_3}^{\text{gen},(n)}(Z) = \sum_{\{K\} \in \text{gen}(M_1,M_2,M_3)} \frac{\Theta^{(n)}(K,Z)}{|O(K)|}$$

where the summation is over a set of representatives of the classes in $\text{gen}(M_1, M_2, M_3)$, $O(K)$ is the (finite) group of orthogonal units of the quadratic lattice $K$, $Z$ is a variable in the Siegel upper half space $\mathbf{H_n}$ and

$$\Theta^{(n)}(K,Z) = \sum_{\boldsymbol{x}=(x_1,\dots,x_n)\in K^n} \exp(2\pi i \text{tr}(q(\boldsymbol{x})Z))$$

with $q(\boldsymbol{x}) = (\frac{1}{2}(M_3 \text{tr}(x_i \bar{x_j})))_{i,j}$.

We consider a double coset decomposition

$$D_{\mathbf{A}}^{\times} = \bigcup_{i=1}^{h=h(M_1,M_2)} D_{\mathbf{Q}}^{\times} y_i R_{\mathbf{A}}^{\times}$$

of the adelic multiplicative group of $D$ with $R_{\mathbf{A}}^{\times} = D_{\infty}^{\times} \times \prod_{p \neq \infty} R_p^{\times}$ and representatives $y_i$ with $n(y_i) = 1$ and $(y_i)_{\infty} = 1$. Then the lattices $I_{ij} = y_i R y_j^{-1}$ (with the norm form scaled by $M_3$) exhaust $\mathrm{gen}(M_1, M_2, M_3)$ (with some classes possibly occurring more than once) and it is easily seen that with $R_i = I_{ii}$ and $e_i = |R_i^{\times}|$ we have

$$\sum_{i,j=1}^{h} \frac{\Theta^{(n)}(I_{ij}, Z)}{e_i e_j} = 2^{\omega} \Theta_{M_1,M_2,M_3}^{\mathrm{gen},(n)}(Z) \tag{5.1}$$

where $\omega = \omega(M_1, M_2)$ is the number of prime divisors of $M_1 M_2$. With these notations we have from [2] (Theorem 3.2 and Corollary 3.2) and [5], (p.229):

**Lemma 5.1** *The value at $s = 0$ of $G_{2,s}(Z)$ is*

$$\sum_{M_1 M_2 M_3 | N} \alpha_{M_1,M_2,M_3} \Theta_{M_1,M_2,M_3}^{\mathrm{gen},(3)}(Z)$$

*with*

$$\alpha_{M_1,M_2,M_3} = (-1)^{1+\omega(M_1,M_2)} (M_1 M_2)^{-3} M_3^{-6} 8\pi^4 \zeta^{(N)}(2)^{-2}.$$

In order to compute the value at $s = 0$ of the differentiated Eisenstein series from Section 2 we have to compute

$$\left( \mathcal{D}_2^{*(a,\nu_2,\nu_3)} \Theta^{(3)}(K, -) \right)(\iota_{111}(z_1, z_2, z_3)) \tag{5.2}$$

for the individual theta series appearing in the sum in Lemma 5.1. We denote by $U_{\mu}$ the space of homogenous harmonic polynomials of degree $\mu$ in 4 variables and identify an element of $U_{\mu}$ with a polynomial on $D_{\infty}$ by evaluating it at the component vector of an element of $D_{\infty}$ with respect to an orthonormal basis relative to the quaternion norm on $D$. Similarly, for $\nu \in \mathbf{N}$ let $U_{\nu}^{(0)}$ be the space of homogeneous harmonic polynomials of degree $\nu$ on $\mathbf{R}^3$ and view $P \in U_{\nu}^{(0)}$ as a polynomial on $D_{\infty}^{(0)} = \{x \in D_{\infty} | tr(x) = 0\}$. The representations $\tau_{\nu}$ of $D_{\infty}^{\times}/\mathbf{R}^{\times}$ of highest weight $(\nu)$ on $U_{\nu}^{(0)}$ given by $(\tau_{\nu}(y))(P)(x) = P(y^{-1}xy)$ for $\nu \in \mathbf{N}$ give all the isomorphism classes of irreducible rational representations of $D_{\infty}^{\times}/\mathbf{R}^{\times}$. By $\langle\langle \quad , \quad \rangle\rangle_{\nu}^{(0)}$ we denote the invariant scalar product in the representation space $U_{\nu}^{(0)}$, by $\langle\langle \quad , \quad \rangle\rangle_{\mu}$ the invariant scalar product in the $SO(D_{\infty}, \mathrm{norm}) =: H_{\mathbf{R}}^{+}$-space $U_{\mu}$. We notice that the invariant scalar products $\langle\langle , \rangle\rangle_{\mu_i}$ on the $U_{\mu_i}$ can be normalized in such a way that they take rational values on the subspaces of polynomials with rational coefficients and that these subspaces generate the $U_{\mu_i}$. Indeed, consider the Gegenbauer polynomial $C^{(\mu_i)}(x, x') =$ obtained from

$$C_1^{(\mu_i)}(t) = 2^{\mu_i} \sum_{j=0}^{[\frac{\mu_i}{2}]} (-1)^j \frac{1}{j!(\mu_i - 2j)!} \frac{(\mu_i - j)!}{2^{2j}} t^{\mu_i - 2j}$$

by

$$\tilde{C}^{(\mu_i)}(x_1, x_2) = 2_i^{\mu}(\mathrm{norm}(x_1)\mathrm{norm}(x_2))^{\mu_i/2}C_1^{(\mu_i)}(\frac{\mathrm{tr}(x_1\overline{x_2})}{2\sqrt{\mathrm{norm}(x_1)\mathrm{norm}(x_2)}})$$

and normalize the scalar product on $U_{\mu_i}$ such that $C^{(\mu_i)}$ is a reproducing kernel, i. e.

$$\langle\langle C^{(\mu_i)}(x, x'), Q(x)\rangle\rangle_{\mu_i} = Q(x')$$

for all $Q \in U_{\mu_i}$. Then the $C^{(\mu_i)}(\cdot, x')$ with rational $x'$ are rational, generate $U_{\mu_i}$, and the reproducing property implies that they have rational scalar products whith each other. The same argument applies to the $U_\nu^{(0)}$.

It is well known that the group of proper similitudes of the quadratic space $(D, \mathrm{norm})$ is isomorphic to $(D^\times \times D^\times)/Z(D^\times)$ via

$$(x_1, x_2) \mapsto \sigma_{x_1, x_2} \text{ with } \sigma_{x_1, x_2}(y) = x_1 y x_2^{-1}$$

and that under this isomorphism $SO(D, \mathrm{norm})$ is the image of

$$\{(x_1, x_2) \in D^\times \times D^\times \mid n(x_1) = n(x_2)\}.$$

Moreover, the $SO(D_\infty, \mathrm{norm}) =: H_{\mathbf{R}}^+$-space $U_\nu^{(0)} \otimes U_\nu^{(0)}$ is isomorphic to the $H_{\mathbf{R}}^+$-space $U_{2\nu}$ and the isomorphism can be normalized in such a way that it preserves rationality and is compatible with the invariant scalar products on both spaces (which are assumed to be normalized as above).

Denoting by $S$ the Gram matrix of the quadratic lattice $K$ we know from Section 1 that form (5.2) is of the form

$$\sum_{\boldsymbol{x}=(x_1, x_2, x_3)\in(\mathbf{Z}^4)^3} P(S^{\frac{1}{2}}x_1, S^{\frac{1}{2}}x_2, S^{\frac{1}{2}}x_3) \exp(\pi i (S[x_1]z_1 + S[x_2]z_2 + S[x_3]z_3))$$

where $P \in \otimes_{i=1}^3 U_{k_i-2}$ is a harmonic polynomial of degree $\mu_i = k_i - 2$ on $\mathbf{R^4}$ in each of the variables and is invariant under the (diagonal) action of $H_{\mathbf{R}} = O(D \otimes \mathbf{R}, \mathrm{norm})$. Moreover, $P$ is independent of $S$ and has rational coefficients (up to a factor of $\pi^{3a+2b}$). The $H_{\mathbf{R}}$-invariant trilinear form $T$ on $U_{\mu_1} \otimes U_{\mu_2} \otimes U_{\mu_3}$ defined by taking the scalar product with the invariant polynomial $\pi^{-3a-2b}P$ ( as remarked in Section 1 this is up to scalars the unique invariant trilinear form) is hence rational (i. e. takes rational values on tensor products of polynomials with rational coefficients).

If all the $\mu_i$ are even (as is the case in our situation) then the decomposition of $H_{\mathbf{R}}^+$ and of $U_\mu$ as $U_{\mu/2}^{(0)} \otimes U_{\mu/2}^{(0)}$ from above gives furthermore that $T$ factors as $T^{(0)} \otimes T^{(0)}$, where the unique (up to scalars) $D_\infty^\times$-invariant trilinear form $T^{(0)}$ on $U_{\mu_1/2}^{(0)} \otimes U_{\mu_2/2}^{(0)} \otimes U_{\mu_3/2}^{(0)}$ has the same rationality properties as $T$. Of course both $T$ and $T_0$ are just the ordinary multiplication if all the $\mu_i$ are 0.

For any positive definite symmetric $4 \times 4$-matrix $S$ we define the $U_{\mu_i}$-valued theta series $\tilde{\Theta}_S^{(\mu_i)}$ by

$$\tilde{\Theta}_S^{(\mu_i)}(z)(x') = \sum_{x\in\mathbf{Z}^4} C^{(\mu_i)}(S^{1/2}x, S^{1/2}x') \exp(\pi i S[x]z) \tag{5.3}$$

We notice that if $K$ is a quaternary quadratic lattice with Gram matrix $S$ the right hand side of 5.3 does not depend on the choice of basis of $K$ with respect to which the Gram matrix is computed (because of the invariance of $C^{(\mu_i)}$ under the (diagonal) action of the orthogonal group); we may therefore write it as $\tilde{\Theta}^{(\mu_i)}(K)$ as well.

We denote by $\tilde{\Theta}^{(\mu_i)}_{M_1,M_2,M_3}$ the weighted average of the $\tilde{\Theta}^{(\mu_i)}_S$ over the Gram matrices of representatives of the classes in the genus gen$(M_1, M_2, M_3)$ as above. We find

$$(\mathcal{D}_2^{*(a,\nu_2,\nu_3)}\Theta(S,-))(\iota_{111}(z_1,z_2,z_3)) = \pi^{3a+2b}T(\tilde{\Theta}^{(\mu_1)}(S,z_1) \otimes \tilde{\Theta}^{(\mu_2)}(S,z_2) \otimes \tilde{\Theta}^{(\mu_3)}(S,z_3))$$
$$(5.4)$$

(where $T$ is as above). For the value at $s = 0$ of (2.1) we obtain therefore

$$\pi^{3a+2b}\sum_{M_1M_2M_3|N}\alpha_{M_1,M_2,M_3}$$
$$\times T(\langle\langle\tilde{\Theta}^{(\mu_1)}_{M_1,M_2,M_3}(z_1),f(z_1)\rangle \otimes \langle\tilde{\Theta}^{(\mu_2)}_{M_1,M_2,M_3}(z_2),\phi(z_2)\rangle \otimes \langle\tilde{\Theta}^{(\mu_3)}_{M_1,M_2,M_3}(z_3),\psi(z_3)\rangle)) \quad (5.5)$$

where by $\langle,\rangle$ we denote the Petersson product.

In order to evaluate this expression further we use Eichler's correspondence. We fix $M_1$, $M_2$ and an Eichler order $R(M_1,M_2) \subseteq D(M_1)$ as above and set $M = M_1M_2$. For an irreducible rational representation $(V_\tau,\tau)$ (with $\tau = \tau_\nu$ as above) of $D^\times_{\mathbf{R}}/\mathbf{R}^\times$ we denote by $\mathcal{A}(D^\times_{\mathbf{A}},R^\times_{\mathbf{A}},\tau)$ the space of functions $\varphi : D^\times_{\mathbf{A}} \to V_\tau$ satisfying $\varphi(\gamma x u) = \tau(u_\infty^{-1})\varphi(x)$ for $\gamma \in D^\times_{\mathbf{Q}}$ and $u = u_\infty u_f \in R^\times_{\mathbf{A}}$. It has been discovered by Eichler that these functions are in correspondence with the elliptic modular forms of weight $2 + 2\nu$ and level $M = M_1M_2$. This correspondence can be described as follows (using the double coset decomposition $D^\times_{\mathbf{A}} = \bigcup_{i=1}^{h} D^\times y_i R^\times_{\mathbf{A}}$ from above): Recall from Section 5 of [2] and Section 3 of [4] that for each $p \mid M_1M_2$ we have an involution $\widetilde{w_p}$ on the space $\mathcal{A}(D^\times_{\mathbf{A}},R^\times_{\mathbf{A}},\tau_\nu)$ given by right translation by a suitable element $\pi_p \in R^\times_p$ of norm $p$ normalizing $R_p$. This space then splits into common eigenspaces of all these (pairwise commuting) involutions. On $\mathcal{A}(D^\times_{\mathbf{A}},R^\times_{\mathbf{A}},\tau_\nu)$ we have furthermore for $p \nmid N$ Hecke operators $\tilde{T}(p)$ whose action on these functions is expressed by the Brandt-matrices $(B^{(\nu)}_{ij}(p))$ (whose entries are endomorphisms of $V_{\tau_\nu}$). They commute with the involutions $\widetilde{w_p}$. On the space $\mathcal{A}(D^\times_{\mathbf{A}},R^\times_{\mathbf{A}},\tau_\nu)$ we have moreover the natural inner product $\langle \ , \ \rangle_\nu$ defined by integration, it is explicitly given by

$$\langle\varphi,\rho\rangle_\nu = \sum_{i=1}^{h}\frac{\langle\langle\varphi(y_i),\rho(y_i)\rangle\rangle_\nu^{(0)}}{e_i}.$$

By abuse of language we call (in the case $\nu = 0$) forms cuspidal, if they are orthogonal to the constant functions with respect to this inner product.

We denote for $p$ dividing $M_2$ the $p$-essential part by $\mathcal{A}_{p,\text{ess}}(D^\times_{\mathbf{A}},R^\times_{\mathbf{A}},\tau)$ consisting of functions $\varphi$ that are orthogonal to all $\rho \in \mathcal{A}(D^\times_{\mathbf{A}},(R'_{\mathbf{A}})^\times,\tau)$ for orders $R' \supseteq R$ for which the completion $R'_p$ strictly contains $R_p$. It is invariant under the $\tilde{T}(p)$ for $p \nmid M_1M_2$ and the $\tilde{w}_p$ for $p \mid M_1M_2$ and hence has a basis of common eigenfunctions of all the $\tilde{T}(p)$ for $p \nmid N$ and all the involutions $\tilde{w}_p$ for $p \mid N$. Moreover by the results of [7, 15, 23, 17] we know that in the space $\mathcal{A}_{\mathrm{ess}}(D_\mathbf{A}^\times, R_\mathbf{A}^\times, \tau)$ of forms that are $p$-essential for all $p$ dividing $M = M_1 M_2$ strong multiplicity one holds, i.e., each system of eigenvalues of the $\tilde{T}(p)$ for $p \nmid M$ occurs at most once, and the eigenfunctions are in one to one correspondence with the newforms in the space $S^{2+2\nu}(M)$ of elliptic cusp forms of weight $2 + 2\nu$ for the group $\Gamma_0(M)$ that are eigenfunctions of all Hecke operators (if $\tau$ is the trivial representation and $R$ is a maximal order one has to restrict here to cuspidal forms on the quaternion side in order to obtain cusp forms on the modular forms side). This correspondence (Eichler's correspondence) preserves Hecke eigenvalues for $p \nmid M$, and if $\varphi$ corresponds to $g \in S^{2+2\nu}(M)$ then the eigenvalue of $g$ under the Atkin-Lehner involution $w_p$ is equal to that of $\varphi$ under $\tilde{w}_p$ if $D$ splits at $p$ and equal to minus that of $\varphi$ under $\tilde{w}_p$ if $D_p$ is a skew field. From (3.13) of [5] we know that if $g$ having first Fourier coefficient 1 corresponds in this way to $\varphi$ with $\langle \varphi, \varphi \rangle_\nu = 1$ then $\langle g, \tilde{\Theta}^{(\mu)}(K) \rangle = \langle g, g \rangle (\varphi(y_i) \otimes \varphi(y_j))$ holds.

It is not difficult (see also [13]) to extend this correspondence to not necessarily new forms $g \in S^{2+2\nu}(N)$ in the following way:

**Lemma 5.2** *Let $N = M_1 M_2$ be a decomposition as above. Call $\tilde{g} \in S^{2+2\nu}(N)$ an $M'$-new form if it is orthogonal to all oldforms coming from $g' \in S^{2+2\nu}(M')$ for $M' \mid N$. Then Eichler's correspondence from above extends to a one to one correspondence between the set of all $M_1$-new eigenforms of the Hecke operators for $p \nmid N$ in $S^{2+2\nu}(N)$ that are eigenfunctions of all the $w_p$ for $p \mid M_2$ with the set of all Hecke eigenfunctions in $\mathcal{A}_t(D_\mathbf{A}^\times, R_\mathbf{A}^\times, \tau)$ that are eigenfunctions of all the involutions $\tilde{w}_p$. This correspondence is compatible with the Hecke action and the eigenvalues under the respective involutions as above; it maps newforms of level $M' \mid N$ (with $M_1 \mid M'$) to forms that are $p$-essential precisely for the $p \mid (M'/M_1)$. The correspondence can be explicitly given (in a nonlinear way) by the first Yoshida lifting sending $\tilde{\varphi}$ to the form*

$$\sum_{ij}^h \frac{\langle \langle \tilde{\varphi}(y_i) \otimes \tilde{\varphi}(y_j), \tilde{\Theta}^{(2\nu)}(I_{ij}) \rangle \rangle}{e_i e_j},$$

*it then sends $\tilde{\varphi}$ with $\langle \tilde{\varphi}, \tilde{\varphi} \rangle_\nu = 1$ to $\tilde{g}$ having first Fourier coefficient one. Moreover, in this normalization it satisfies the scalar product relations from above, i. e. , if $\tilde{g}$ corresponds to $\tilde{\varphi}$ then $\langle \tilde{g}, \tilde{\Theta}^{(2\nu)}(I_{ij}) \rangle \langle \tilde{\varphi}, \tilde{\varphi} \rangle_\nu^2 = \langle \tilde{g}, \tilde{g} \rangle (\tilde{\varphi}(y_i) \otimes \tilde{\varphi}(y_j))$ holds.*

*Proof.* Let $M'$ be a divisor of $N$ and $g \in S^{2+2\nu}(M')$ and let $\epsilon$ be a function from the set $S$ of prime divisors of $N/M'$ (whose cardinality we denote by $\omega(M/M')$) to $\{\pm 1\}^{\omega(M/M')}$. Then the function $g^\epsilon := \sum_{S' \subseteq S} g | \prod_{p \in S} \epsilon(p) w_p$ in $S^{2+2\nu}(N)$ is an eigenfunction of the Atkin-Lehner involutions $w_p$ for the $p \mid (N/M')$ with eigenvalues $\epsilon(p)$. Similarly, fix a maximal order $\tilde{R} \subseteq D$ and an Eichler order $R = R(M_1, M_2) \subseteq \tilde{R}$, let $M_2' \mid M_2$ and let $R(M_1, M_2')$ be the Eichler order of level $M_1 M_2'$ in $D$ containing $R$ and contained in $\tilde{R}$. Let $\epsilon$ be a function on the set of prime divisors of $M_2/M_2'$ as above. Then to $\varphi \in \mathcal{A}_{\mathrm{ess}}(D_\mathbf{A}^\times, (R(M_1, M_2'))_\mathbf{A}^\times, \tau)$ we construct as above a unique $\varphi^\epsilon$ having the same Hecke eigenvalues for $p \nmid N$ and the same $\tilde{w}_p$-eigenvalues for $p \mid M_1 M_2'$ as $\varphi$ such that $\varphi^\epsilon$ is an eigenfunction of the $\tilde{w}_p$ for $p \mid (M_2/M_2')$ with eigenvalues $\epsilon(p)$.

Given an $M_1$-new Hecke eigenform $\tilde{g}$ in $S^{2+2\nu}(N)$ that is an eigenfunction of all the $w_p$ for $p \mid M_2$ with eigenvalues $\epsilon_p$ we then associate to it the newform $g$ of some level $M' = M_1 M_2' \mid M_1 M_2$ that has the same Hecke eigenvalues for $p \nmid N$ so that $\tilde{g} = g^\epsilon$ with $\epsilon(p) = \epsilon_p$ and apply Eichler's correspondence to get a $\varphi \in \mathcal{A}_{\mathrm{ess}}(D_\mathbf{A}^\times, (R(M_1, M_2'))_\mathbf{A}^\times, \tau_\nu)$. From [2] we know that this can be normalized such that $g$ is obtained from $\varphi$ by Yoshida's lifting. We then pass to the eigenfunction $\varphi^\epsilon$ of all $\widetilde{w_p}$ with the same eigenvalues as $g$ for $p \mid M_2$. The scalar product relation (up to normalization) follows then in the same way as in [5], using the uniqueness of the given set of $\tilde{w}_p$-eigenvalues and the fact that application of $w_p$ for $p \mid M_2$ transforms $\tilde{\Theta}^{(2\nu)}(I_{ij})$ to $\tilde{\Theta}^{(2\nu)}(I_{i'j})$, where $y_{i'}$ represents the double coset of $y_i \pi_p^{-1}$ (this is an easy generalization of Lemma 9.1 a) of [2], see also [3]). The same argument shows that Yoshida's lifting realizes this correspondence (using the well known fact that it gives the right Hecke eigenvalues for the $p \nmid N$), and using the expression of $g^\epsilon$ as a Yoshida lifting we find the correct normalization of the scalar product relation as in [5].

We will need a version of the scalar product relation in Lemma 5.2 also for the case of newforms of level strictly dividing the level of the theta series involved.

**Lemma 5.3** *Let $M_1, M_2'$ with $M = M_1 M_2'$ dividing $N$ be as before and let $g$ be a normalized newform of level $M$ and weight $k = 2 + 2\nu$ as in the previous Lemma; let $S$ be the set of prime divisors of $N/M$. Put $M_2 = N/M_2'$ and let $R' = R(M_1, M_2')$ and $R = R(M_1, M_2) \subseteq R'$ be Eichler orders of levels $M, N$ respectively in $D(M_1)$ and consider a double coset decomposition $D_\mathbf{A}^\times = \cup_{i=1}^h D_\mathbf{Q}^\times y_i R_\mathbf{A}^\times$ and corresponding quadratic lattices (ideals in $D$) $I_{ij}$ relative to $R$ as before. Let $\varphi \in \mathcal{A}_{\mathrm{ess}}(D_\mathbf{A}^\times, (R'_\mathbf{A})^\times, \tau_\nu)$ be the essential form corresponding to $g$ under Eichler's correspondence (with $\langle \varphi, \varphi \rangle_\nu = 1$). Then*

$$\langle g, \tilde{\Theta}^{(2\nu)}(I_{ij}) \rangle = \langle g, g \rangle \sum_{S' \subseteq S} \Big( \prod_{p \in S'} \widetilde{(w_p)} \varphi \Big)(y_i) \otimes \Big( \prod_{p \in S'} \widetilde{(w_p)} \varphi \Big)(y_j) \tag{5.6}$$

*Proof.* Let $\epsilon$ be a function from the set of prime divisors of $N/M$ to $\{\pm 1\}^{\omega(N/M)}$ and let $g^\epsilon, \varphi^\epsilon$ be as in the proof of Lemma 5.2. Then $g^\epsilon$ corresponds to $(\varphi^\epsilon)/\sqrt{\langle \varphi^\epsilon, \varphi^\epsilon \rangle_\nu}$ under Yoshida's correspondence, so we get

$$\langle g^\epsilon, \tilde{\Theta}^{(2\nu)}(I_{ij}) \rangle = \frac{\langle g^\epsilon, g^\epsilon \rangle}{\langle \varphi^\epsilon, \varphi^\epsilon \rangle_\nu} \varphi^\epsilon(y_i) \otimes \varphi^\epsilon(y_j) \tag{5.7}$$

by Lemma 5.2. For any $S' \subseteq S$ the Petersson product $\langle g, g | (\prod_{p \in S'} \epsilon(p)(w_p)) \rangle$ is the same as the product of $g$ with the image of $g | (\prod_{p \in S'} \epsilon(p)(w_p))$ under the trace operator from modular forms for $\Gamma_0(N)$ to modular forms for $\Gamma_0(M)$, hence equal to $\prod_{p \in S'} \epsilon(p) p^{-\nu} a_g(p) \langle g, g \rangle$, since for $p \nmid M$ the map $g \mapsto g | w_p$ composed with the trace just gives the (renormalized) Hecke operator. The same argument applies to $\langle \varphi, (\prod_{p \in S'} \epsilon(p) \widetilde{(w_p)} \varphi) \rangle$ and gives the same factor of comparison with $\langle \varphi, \varphi \rangle$ since $g$ and $\varphi$ have the same (renormalized) Hecke eigenvalues. Thus we have $\langle g^\epsilon, g^\epsilon \rangle \langle \varphi^\epsilon, \varphi^\epsilon \rangle^{-1} = \langle g, g \rangle \langle \varphi, \varphi \rangle^{-1}$, and summing up the identities (5.7) for all the functions $\epsilon$ gives the assertion.

**Lemma 5.4** *Let $f, \phi, \psi$ as in Section 2 be newforms of square free levels $N_f, N_\phi, N_\psi$ with $N = lcm(N_f, N_\phi, N_\psi)$ and with weights $k_i = 2 + 2\mu_i$. Then in 5.5 the summand for $M_1, M_2, M_3$ is zero unless $M_3 = 1$, $N = M_1 M_2$ and $M_1 \mid \gcd(N_f, N_\phi, N_\psi)$ hold.*

*Proof.* Since by Lemma 5.2 we can express $f, \phi, \psi$ as Yoshida-liftings an easy generalization of Lemma 9.1 b) of [2] shows that $f$ is orthogonal to all $\tilde{\Theta}^{(\mu_1)}(K)$ for $K$ in $\text{gen}(M_1, M_2, M_3)$ for which $N_f$ does not divide $M_1 M_2$ and analogously for $\phi, \psi$. This establishes the vanishing of all summands for which $N \neq M_1 M_2$.

If there is a $p \mid M_1$ not dividing $\gcd(N_f, N_\phi, N_\psi)$ then say $p \nmid N_f$. The Petersson product of $f$ with $\tilde{\Theta}^{(2\mu_1)}(K)$ for $K$ in $\text{gen}(M_1, M_2, 1)$ is then the same as that of $f$ with the form obtained by applying the trace operator from modular forms on $\Gamma_0(N)$ to modular forms on $\Gamma_0(N/p)$ to the theta series. But it is easily checked that this trace operator annihilates the theta series of $K \in \text{gen}(M_1, M_2, M_3)$ if $p \mid M_1$, see [8]; the same argument is applied to $\phi, \psi$ which shows the last part of the assertion.

**Lemma 5.5** *Let $\tilde{f}, \tilde{\phi}, \tilde{\psi}$ be cusp forms of weights $k_1, k_2, k_3$ for $\Gamma_0(N)$ with square free $N$ as in Section 2, assume them to be eigenfunctions of the Hecke operators for $p \nmid N$ and of all the Atkin-Lehner involutions $w_p$ with eigenvalues $\epsilon_f(p), \epsilon_\phi(p), \epsilon_\psi(p)$ but not necessarily newforms. Let $M_1 M_2 = N$ (with $M_1$ as always having an odd number of prime factors) and let $R = R(M_1, M_2)$ be an Eichler order in $D = D(M_1)$. Let $\varphi_f, \varphi_\phi, \varphi_\psi$ be the forms in $\mathcal{A}(D_\mathbf{A}^\times, R_\mathbf{A}^\times, \tau_{\nu_i})$ corresponding to $f, \phi, \psi$ under the correspondence of Lemma 5.2 (with $k_i = 2 + 2\mu_i$ for $i = 1, \ldots, 3$). Then the summand for $M_1, M_2, M_3 = 1$ in 5.5 is*

$$2^{-\omega(N)} \alpha_{M_1, M_2, 1} \langle f, f \rangle \langle \phi, \phi \rangle \langle \psi, \psi \rangle \sum_{i=1}^{h} \left( \frac{T_0(\varphi_f(y_i) \otimes \varphi_\phi(y_i) \otimes \varphi_\psi(y_i))}{e_i} \right)^2.$$

*The latter expression is zero unless for all $p \mid N$ one has $p \mid M_1$ if and only if $\epsilon_f(p)\epsilon_\phi(p)\epsilon_\psi(p) = -1$.*

*Proof.* The first part of the assertion is an immediate consequence of Lemma 5.2 and the decomposition of $T$ as $T^{(0)} \otimes T^{(0)}$. For the second part we notice that the expression $\sum_{i=1}^{h} \left( \frac{T_0(\varphi_f(y_i) \otimes \varphi_\phi(y_i) \otimes \varphi_\psi(y_i))}{e_i} \right)$ does not change if an involution $\widetilde{w_p}$ is applied to all three functions $\varphi_f, \varphi_\phi, \varphi_\psi$ since this only permutes the order of summation. On the other hand each summand is multiplied with the product of the eigenvalues of $\varphi_f, \varphi_\phi, \varphi_\psi$ under $\widetilde{w_p}$, which in view of the relation between the $w_p$ eigenvalues and the $\widetilde{w_p}$-eigenvalues of corresponding functions proves the assertion.

Although the statement of Lemma 5.2 is not true if one omits the condition that $\tilde{f}$ is an eigenfunction of all the involutions $w_p$ the next Lemma shows that by an amusing newforms argument an only slightly changed version of Lemma 5.5 remains true without this condition.

**Lemma 5.6** *Let $f, \phi, \psi$ be normalized newforms of levels $N_f, N_\phi, N_\psi$ of weights $k_i = 2 + 2\mu_i$ ($i = 1, \ldots, 3$) as in Section 2 and let $\Lambda_1, \Lambda_2, \Lambda_3$ be pairwise disjoint subsets of the set of primes divisors of $N^f = N/N_f, N^\phi, N^\psi$ respectively such that $\Lambda_1 \cup \Lambda_2 \cup \Lambda_3$ is the set of all primes* *dividing precisely one of the integers* $N_f, N_\phi, N_\psi$. *For* $\kappa = 1, \ldots, 3$ *let*

$$\widetilde{w_{\Lambda_\kappa}} = \prod_{p \in \Lambda_\kappa} \widetilde{w_p}.$$

*Let a decomposition* $N = M_1 M_2$ *as before be given, fix a maximal order* $\tilde{R}$ *in* $D = D(M_1)$ *and an Eichler order* $R = R(M_1, M_2) \subseteq \tilde{R}$ *of level* $M_1 M_2 = N$ *in* $D$ *and consider a double coset decomposition* $D_\mathbf{A}^\times = \cup_{i=1}^h D_\mathbf{Q}^\times y_i R_\mathbf{A}^\times$ *and corresponding quadratic lattices (ideals in* $D$) $I_{ij}$ *relative to* $R$ *as before. For each* $M_2' \mid M_2$ *let* $R(M_1 M_2') = R(M_1, M_2')$ *be the unique Eichler order of level* $M_1 M_2'$ *contained in* $\tilde{R}$ *and containing* $R$. *Let* $\varphi_1 \in \mathcal{A}_{ess}(D_\mathbf{A}^\times, (R(N_f)_\mathbf{A}^\times, \tau_{\mu_1}))$ *be the form corresponding to* $f$ *under Eichler's correspondence and define* $\varphi_2, \varphi_3$ *analogously with respect to* $\phi, \psi$. *Then*

$$T\Big(\langle f, \tilde{\Theta}^{(\mu_1)} \rangle \langle \phi, \tilde{\Theta}^{(\mu_2)} \rangle \langle \psi, \tilde{\Theta}^{(\mu_3)} \rangle\Big)$$
$$= 2^{-\omega(\gcd(N_f, N_\phi, N_\psi))} \langle f, f \rangle \langle \phi, \phi \rangle \langle \psi, \psi \rangle \Big(T_0\Big(\sum_{i=1}^h \frac{1}{e_i} \widetilde{w_{\Lambda_1}} \varphi_1(y_i) \otimes \widetilde{w_{\Lambda_2}} \varphi_2(y_i) \otimes \widetilde{w_{\Lambda_3}} \varphi_3(y_i)\Big)\Big)^2 \quad (5.8)$$

*Proof.* This is an immediate consequence of Lemmas 5.2 and 5.3: Upon inserting the scalar product relations from these lemmata into the left hand side of (5.8) we obtain a sum of terms of the type

$$\Big(T_0\big((\sum_{i=1}^h \widetilde{w_{\Lambda_1'}} \varphi_1(y_i) \otimes \widetilde{w_{\Lambda_2'}} \varphi_2(y_i) \otimes \widetilde{w_{\Lambda_3'}} \varphi_3(y_i)\big)\Big)^2$$

with arbitrary subsets $\Lambda_\kappa'$ of the sets of primes dividing $N/N_f, N/N_\phi, N/N_\psi)$ respectively. Let $p$ be a prime dividing two of the levels, say $p \mid N_f, p \mid N_\phi$. Then since applying $\widetilde{w_p}$ to all three of the $\varphi_\kappa$ only changes the order of summation, the involution $\widetilde{w_p}$ for $p \in \Lambda_3'$ may be pulled over to $\varphi_2, \varphi_3$ which are eigenfunctions of $\widetilde{w_p}$. The terms with the set $\Lambda_3'$ and those with $\Lambda_3' \setminus \{p\}$ give therefore the same contribution. Let now $p$ be a prime dividing only one of the levels, say $p \mid N_f$. If $p \notin \Lambda_2' \cup \Lambda_3'$ then the component of $\varphi_2 \otimes \varphi_3$ in the $\tau_{\mu_1}$-isotypic component of $\tau_{\nu_2} \otimes \tau_{\nu_3}$ is an oldform with respect to $p$, hence orthogonal to the $p$-essential form $\varphi_1$. Since $T_0(\varphi_1 \otimes \varphi_2 \otimes \varphi_3)$ is proportional to the scalar product of $\varphi_1$ with this component of $\varphi_1 \otimes \varphi_2$ such a term gives no contribution; the same argument applies if $p \in \Lambda_2' \cap \Lambda_3'$ holds. If $p$ is in precisely one of $\Lambda_2', \Lambda_3'$ then the same argument as in the first case shows that $p$ may be shifted to either one of these sets without changing the contribution of the term. Taking together both cases we find that all terms appearing are of the shape

$$\Big(T_0\big((\sum_{i=1}^h \widetilde{w_{\Lambda_1}} \varphi_1(y_i) \otimes \widetilde{w_{\Lambda_2}} \varphi_2(y_i) \otimes \widetilde{w_{\Lambda_3}} \varphi_3(y_i)\big)\Big)^2$$

with each term appearing $2^{\omega(N/\gcd(N_f, N_\phi, N_\psi))}$ times, which in view of 5.1 implies the assertion.

Collecting all the information obtained we arrive at the main theorem:

**Theorem 5.7** *The value of the triple product $L$-function $L(f, \phi, \psi, s)$ at the central critical value $s = \frac{k_1 + k_2 + k_3}{2} - 1$ is*

$$(-1)^{a'} 2^{5+4a+3b-\omega(\gcd(N_f, N_\phi, N_\psi))} \pi^{5+9a'+4b} \frac{(a'+1)^{[b]}}{2^{[a+b]} 2^{[a']} (\nu_2 + 1)^{[a']} (\nu_3 + 1)^{[a']}}$$

$$\times \langle f, f \rangle \langle \phi, \phi \rangle \langle \psi, \psi \rangle \Big( T_0 \Big( \sum_{i=1}^{h} \frac{1}{e_i} \widetilde{w_{\Lambda_1}} \varphi_1(y_i) \otimes \widetilde{w_{\Lambda_2}} \varphi_2(y_i) \otimes \widetilde{w_{\Lambda_3}} \varphi_3(y_i) \Big) \Big)^2 \quad (5.9)$$

*where the notation is as in Lemma 5.6 and $T_0$ is (as explained in the beginning of this section) the up to scalars unique rational invariant trilinear form on the representation space $U_{\mu_1/2}^{(0)} \otimes U_{\mu_2/2}^{(0)} \otimes U_{\mu_3/2}^{(0)}$ (with $\mu_i = k_i - 2$) and takes values in the coefficient fields of $f, \phi, \psi$ respectively on the polynomials $\widetilde{w_{\Lambda_1}} \varphi_\kappa(y_i)$ (for $\kappa = 1, \dots, 3$)*

It should be noted that the rational quantity on the right hand side can be interpreted as the height pairing of a diagonal cycle with itself in the same way as in [11]. One has just to replace (for $\kappa = 1, \dots, 3$) the group $\mathrm{Pic}(X)$ of [11] with the group $\mathrm{Pic}(V_\kappa)$ from [14] obtained by attaching to each $y_i$ in the double coset decomposition of $D_{\mathbf{A}}^{\times}$ used above the space of $R_i^{\times}$-invariant polynomials in $U_{\mu_\kappa}^{(0)}$. Our functions $\varphi_\kappa$ may then be interpreted as elements of $\mathrm{Pic}(V_\kappa)$. One may then form the tensor product of these three groups and obtain an analogue of the diagonal cycle $\Delta$ from [11] by using our Gegenbauer polynomials from above and proceed as in loc. cit.

## References

[1] S.Böcherer: Über die Fourier-Jacobientwicklung der Siegelschen Eisensteinreihen II. Mathem.Z.189(1989),81-110

[2] S.Böcherer,R.Schulze-Pillot: Siegel modular forms and theta series attached to quaternion algebras. Nagoya Math.J. 121(1991), 35-96

[3] S.Böcherer,R.Schulze-Pillot: Siegel modular forms and theta series attached to quaternion algebras II, Preprint 1995

[4] S.Böcherer,R.Schulze-Pillot: Mellin transforms of vector valued theta series attached to quaternion algebras. Math. Nachr. 169 (1994), 31-57

[5] S.Böcherer,R.Schulze-Pillot: Vector valued theta series and Waldspurger's theorem, Abh. Math. Sem. Hamburg 64 (1994), 211-233

[6] S.Böcherer,T.Satoh,T.Yamazaki: On the pullback of a differential operator and its application to vector valued Eisenstein series. Comm.Math.Univ.S.Pauli 41(1992), 1-22

[7] M. Eichler: The basis problem for modular forms and the traces of the Hecke operators, p. 76-151 in Modular functions of one variable *I*, Lecture Notes Math. 320, Berlin-Heidelberg-New York 1973

[8] J. Funke: Spuroperator und Thetareihen quadratischer Formen, Diplomarbeit Köln 1994

[9] P.Garrett: Decomposition of Eisenstein series: Rankin triple products. Annals of Math. 125(1987), 209-235

[10] P. Garrett, M. Harris: Special values of triple product $L$-Functions, Am. J. of Math. 115 (1993), 159-238

[11] B.Gross, S.Kudla: Heights and the central critical values of triple product L-functions.Compositio Math.81(1992), 143-209

[12] M. Harris, S. Kudla: The central critical value of a triple product $L$-function, Annals of Math. 133 (1991), 605-672

[13] K. Hashimoto: On Brandt matrices of Eichler orders, Preprint 1994

[14] R. Hatcher: Heights and *L*-series, Can. J. math. 62 (1990), 533-560

[15] H. Hijikata, H. Saito: On the representability of modular forms by theta series, p. 13-21 in Number Theory, Algebraic Geometry and Commutative Algebra, in honor of Y. Akizuki, Tokyo 1973

[16] T.Ibukiyama: On differential operators on automorphic forms and invariant pluriharmonic polynomials. Preprint 1990

[17] H. Jacquet, R. Langlands: Automorphic forms on $GL(2)$, Lect. Notes in Math. 114, Berlin-Heidelberg-New York 1970

[18] W.Li: Newforms and functional equations. Math.Ann.212(1975), 285-315

[19] P. Littelmann: On spherical double cones, J. of Algebra 166 (1994), 142-157

[20] H.Maaß: Siegel's modular forms and Dirichlet series. Lect. Notes Math. 216 Berlin, Heidelberg, New York: Springer 1971

[21] T.Orloff: Special values and mixed weight triple products (with an appendix by Don Blasius). Invent.math.90(1987), 169-180

[22] T.Satoh: Some remarks on triple L-functions. Math.Ann.276(1987), 687-698

[23] H. Shimizu: Theta series and automorphic forms on $GL_2$. J. of the Math. Soc. of Japan 24 (1972), 638-683

[24] G.Shimura: The special values of the zeta functions associated with cusp forms. Comm.Pure Appl. Math.29(1976), 783-804

[25] G.Shimura: The arithmetic of differential operators on symmetric domains. Duke Math.J.48(1981), 813-843

S. Böcherer, Fakultät für Mathematik und Informatik, Universität Mannheim, Seminargebäude A5, D-68131 Mannheim, Germany

R. Schulze-Pillot, Mathematisches Institut, Universität zu Köln, Weyertal 86-90, D-50931 Köln, Germany
after Oct. 1, 1995: Fachbereich Mathematik, Universität des Saarlandes (Bau 27), Postfach 1150, D-66041 Saarbrücken, Germany